\documentclass[11pt]{article}
	
	\newcommand{\blind}{0}
	
    \makeatletter
    \renewcommand\section{\@startsection {section}{1}{\z@}%
                                       {-3.5ex \@plus -1ex \@minus -.2ex}%
                                       {2.3ex \@plus.2ex}%
                                       {\normalfont\fontfamily{phv}\fontsize{16}{19}\bfseries}}
    \renewcommand\subsection{\@startsection{subsection}{2}{\z@}%
                                         {-3.25ex\@plus -1ex \@minus -.2ex}%
                                         {1.5ex \@plus .2ex}%
                                         {\normalfont\fontfamily{phv}\fontsize{14}{17}\bfseries}}
    \renewcommand\subsubsection{\@startsection{subsubsection}{3}{\z@}%
                                        {-3.25ex\@plus -1ex \@minus -.2ex}%
                                         {1.5ex \@plus .2ex}%
                                         {\normalfont\normalsize\fontfamily{phv}\fontsize{14}{17}\selectfont}}
    \makeatother
	
	\usepackage{amsmath}
	\usepackage{graphicx}
	\usepackage{enumerate}
	\usepackage{xcolor}
	\usepackage{natbib} 
	\usepackage{url} 
	
 \usepackage{booktabs}
\usepackage[ruled]{algorithm2e}
\usepackage{algpseudocode}
\usepackage{tabularx}
\usepackage{multirow}
\usepackage{amssymb}
\usepackage{mathrsfs}
\usepackage{caption}
	
\begin{document}
		
		\def\spacingset#1{\renewcommand{\baselinestretch}%
			{#1}\small\normalsize} \spacingset{1}
		
		\if0\blind
		{
			\title{\bf Joint optimization for production operations considering reworking}
   \author{Yilan Shen$^a$, Boyang Li$^a$, Xi Zhang$^{a*}$\\
        $^a$College of Engineering, Peking University, Beijing, China}
			\date{}
			\maketitle
		} \fi
		
		\if1\blind
		{

            \title{\bf Joint optimization for production operations considering reworking}
      \author{Yilan Shen$^a$, Boyang Li$^a$, Xi Zhang$^{a*}$\\
        $^a$College of Engineering, Peking University, Beijing, China}
			
			\begin{center}
				{\LARGE\bf Joint optimization for production operations considering reworking}
			\end{center}
			\medskip
		} \fi
\begin{abstract}
In pursuit of enhancing the comprehensive efficiency of production systems, our study focused on the joint optimization problem of scheduling and machine maintenance in scenarios where product rework occurs. The primary challenge lies in the interdependence between product \underline{q}uality, machine \underline{r}eliability, and \underline{p}roduction scheduling, compounded by the uncertainties from machine degradation and product quality, which is prevalent in sophisticated manufacturing systems. To address this issue, we investigated the dynamic relationship among these three aspects, named as QRP-co-effect. On this basis, we constructed an optimization model that integrates production scheduling, machine maintenance, and product rework decisions, encompassing the context of stochastic degradation and product quality uncertainties within a mixed-integer programming problem. To effectively solve this problem, we proposed a dual-module solving framework that integrates planning and evaluation for solution improvement via dynamic communication. By analyzing the structural properties of this joint optimization problem, we devised an efficient solving algorithm with an interactive mechanism that leverages \emph{in-situ} condition information regarding the production system's state and computational resources. The proposed methodology has been validated through comparative and ablation experiments. The experimental results demonstrated the significant enhancement of production system efficiency, along with a reduction in machine maintenance costs in scenarios involving rework. \\
\end{abstract}
	\noindent%
	\spacingset{1.0} 

\section{Introduction} \label{s:intro}
In contemporary manufacturing industries, production efficiency and quality are regarded by workshop operators as essential indicators of a production system's operational effectiveness. Substantial research efforts have been dedicated to these areas. In the case of single-variety mass production, adept management of the production system allows operators to achieve a harmonious equilibrium between product quality and production output within a reasonable cost range, employing continuous improvement techniques on the production line or optimizing machine configurations. However, challenges arise when dealing with diverse, sophisticated manufacturing systems like semiconductor, aerospace or heavy machinery production, wherein small batch production of customized products with distinct characteristics predominates, causing deviations from specified quality requirements due to frequent process alterations \citep{mabkhot2020ontology}. Consequently, rectification tasks are often performed on defective components rather than discarding them, primarily driven by economic considerations \citep{flapper2002planning}. For instance, in the wafer fabrication industry, the photolithography process defines the graph of each layer, necessitating stringent quality control measures and subsequent rework for any non-conformance \citep{sha2006dispatching}. Additionally, the stochastic degradation of machines further compounds the impact on product quality characteristics during the manufacturing process. Such occurrences can result in products failing to pass quality inspections, compelling them to be sent back for rework at previous stages of production. This phenomenon is commonplace in sophisticated manufacturing systems, exerting a substantial influence on the original production schedule and ultimately diminishing the overall production throughput of the manufacturing system. Moreover, ensuring the attainment of superior product quality is of paramount importance. Unfortunately, despite the urgency, the first pass yield (FPY) indicator, which gauges the proportion of products meeting quality specifications on the initial attempt, remains suboptimal, resulting in an inability to meet prevailing demands. As a consequence, the need for rework arises, rendering adherence to the original production schedule untenable. In summary, the manufacturing system inherently encompasses stochastic properties, frequently grappling with uncertainties stemming from random disruptive events. The stochastic deterioration of machines profoundly impacts both the timeliness and quality of production outcomes. In addition, unpredictable quality defects introduce further random disruptions that influence system performance and necessitate rescheduling. This poses a challenge when designing system models and conducting performance analyses, particularly for real-time assessments of system performance during production.

\begin{figure*}[t]
\centering
\includegraphics[scale=0.4]{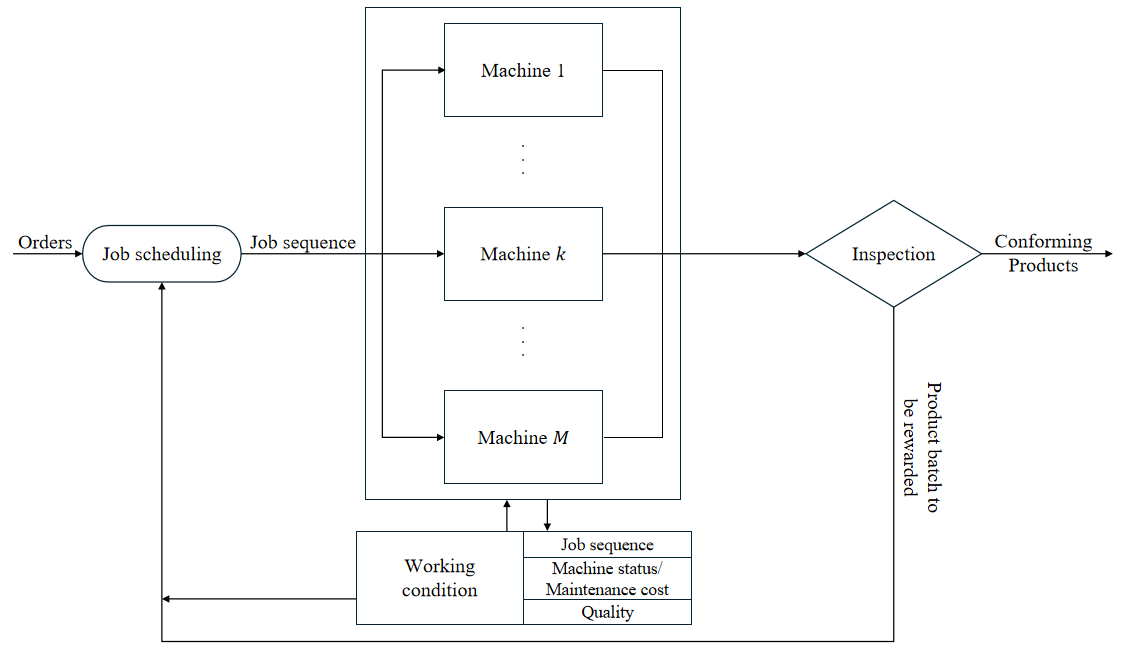}
\caption{\textcolor{black}{The framework of the multi-component system with the reworking activity.}}
\label{Figure_structure}
\end{figure*}

In the realm of production and quality management, the treatment of quality inspection, production capacity, and machine maintenance as separate entities has been the norm. However, current literature has shed light on the profound impact of non-conforming products undergoing rework on the overall performance of production systems \citep{gardner2020managing, amirkhani2017new}. These three elements, in fact, intricately intertwine and exhibit profound dependencies within advanced machine production systems. On one hand, factors such as an increased volume of processing tasks and the presence of non-conforming products resulting from dimensional deviations in tools or processed items can accelerate the degradation rate of each machine. This acceleration subsequently leads to a decline in production rate. On the other hand, machine degradation can result in an upsurge of non-conforming products requiring rework, consequently causing disruptions to the original production planning and task scheduling. For instance, in the context of a single-stage manufacturing process (as depicted in Figure \ref{Figure_structure}), jobs are initially assigned to independent parallel machines based on a predetermined production plan. Nonetheless, during the process, machine deterioration significantly impacts both the actual processing time of jobs and the occurrence of product defects. When the number of inspected non-conforming products exceeds a certain threshold, the jobs necessitating rework are rescheduled, taking into consideration the machine's degradation state and the job sequence to ensure compliance with quality requirements. Furthermore, an appropriate maintenance policy becomes indispensable in addressing the machine's reliability concerns, thereby ensuring prevention of product defects and facilitating a consistent production throughput that yields a relatively high proportion of qualified products. Moreover, the production scheduling itself, as it determines the sequencing of jobs, exerts a profound influence on both the machine's degradation path and the product's quality. Hence, it becomes crucial to meticulously explore the significant interrelationships among quality, production, and maintenance. Such exploration is essential in order to provide a comprehensive strategy for job scheduling and machine maintenance.

Within the extant body of literature, numerous studies have emerged over the past two decades that delve into the intricate relationship among these elements. One pivotal study, originating from \cite{chen2005quality}, pioneered the investigation of the dependency between product quality and machine reliability. This study introduced a model known as the \underline{q}uality and \underline{r}eliability chain (QR chain), which effectively integrated the reliability of manufacturing system components with product quality considerations in multi-stage manufacturing processes. Notably, this work examined the interaction between quality and reliability from the vantage point of system reliability, laying a solid foundation for optimizing machines' operating parameters, such as the locator wear rate, to enhance quality. Building upon this groundbreaking work, subsequent research efforts have devised comprehensive models of integrated analysis, focusing on elucidating the interrelationship between any two of the aforementioned elements. For instance, a representative study by \cite{kang2019integrated} effectively examined the interaction between machine degradation and product throughput. By establishing a correlation between the workload imposed on each machine and the rate at which degradation occurs, this investigation illuminates the potential for optimizing machine degradation by dynamically adjusting workloads. Such optimization ensures the continuity of production throughput within the manufacturing system. Despite these notable contributions, a majority of studies tend to overlook crucial implementation possibilities. Specifically, the stochastic nature of machine degradation paths while processing various jobs from different orders is often neglected. Furthermore, even machines operating under the same workload but with different job sequences exhibit dissimilar degradation paths. This discrepancy in degradation paths necessitates distinct maintenance policies, which, in turn, interfere with job scheduling. Moreover, the simultaneous assurance of both quantity and quality of products assumes paramount importance, thereby necessitating the optimization of production efficiency. In such cases, reworking activities become essential, along with the need to collectively address non-conforming products that require rework while accommodating existing pending jobs. This calls for a more efficacious approach to job sequencing. Unfortunately, prevailing studies tend to propose planning and scheduling solutions for production decisions independently, falling short of providing a comprehensive framework that effectively addresses the uncertain occurrence of rework events stemming from endogenous factors, including machine reliability and production scheduling. 



The advancement of sensor technologies has facilitated the seamless integration of machines into a unified platform, thereby enabling the acquisition of vital \emph{in-situ} process condition information such as job status and product quality. This platform presents an unparalleled opportunity to achieve a paramount objective: the simultaneous determination of maintenance, job arrangement, and reworking activities while considering the intricate interrelationship among quality, production, and machine degradation based on the integration of condition information. By jointly analyzing these three pivotal elements, a more precise schedule can be obtained, a crucial aspect in the realm of sophisticated manufacturing with a primary focus on attaining optimal production efficiency. These compelling factors serve as a significant impetus to undertake a comprehensive study of the interplay among quality, production, and machine degradation within production systems, endeavoring to devise a comprehensive and decisive strategy in terms of maintenance and job scheduling optimization.

While the imperative for quality, delivery efficiency, and reliability remains paramount in the realm of sophisticated manufacturing, there has been a lack of collective investigation into these three elements \citep{cheng2020joint, eben2022relationships} within a production system. The existing literature fails to comprehensively explore the intricate relationship that exists among these components, often oversimplifying their interplay and neglecting the nuances of the job production process \citep{tambe2022reliability}. 

The current understanding of the co-effect of Quality-Reliability-Planning (QRP) emphasizes the pervasive uncertainties surrounding job processing times, machine reliability, and rework activities. These uncertainties are attributed to the stochastic nature of machine degradation, which is influenced by a multitude of factors. This poses a unique challenge in formulating a solution that addresses job scheduling, machine maintenance, and rework, while simultaneously enhancing planning quality and optimizing computational efficiency within a dynamic production system. Adding to the complexity, this solution must effectively account for both endogenous and exogenous uncertainties, providing appropriate production decisions that align harmoniously with real-world production conditions. Unfortunately, the prevailing methodologies lack the necessary adaptability to comprehensively encompass all three elements within uncertain production systems.

To bridge these significant gaps, we develop an integrated approach that takes a holistic perspective, considering all the pivotal elements inherent in sophisticated manufacturing. Our approach ensures the delivery of high-quality products while concurrently maintaining cost-effectiveness. Our contributions are as follows:


1. We propose a novel framework that elucidates the intricate relationship among quality, reliability, and job scheduling. This framework leverages the stochastic nature of sequence-dependent machine degradation, takes into account the uncertain influence of deterioration on job processing time, and considers the impact of job sequence and machine state on product quality characteristics.

2. We specifically focus on the production reworking scenario and establish a model that captures the dynamic changes arising from job scheduling, maintenance, and reworking activities during production.

3. Our solution framework offers flexibility, facilitating production decisions for systems operating under production uncertainties. This framework consists of two interconnected modules: the planning module and the evaluation module. These modules leverage condition information and seamlessly interface through a thoughtfully designed communication mechanism, enabling efficient adjustment of job re-arrangement and machine maintenance.

The ensuing sections of this study are organized as follows: Section 2 presents a succinct review of the relevant literature. Section 3 provides an elaborate description of the problem, incorporating the QRP-co-effect. Section 4 formulates this specific problem based on the QRP-co-effect. Section 5 introduces the proposed two-module solution framework along with the tailored algorithms. In Section 6, an extensive array of numerical studies evaluates the effectiveness of the proposed solution method. Ultimately, Section 7 concludes the study and delves into prospects for future research.




\section{Related literature}

Quality control, production scheduling, and reliability are three fundamental production elements in manufacturing systems \citep{das2007machine, wang1998product, sun2023robust}, which are directly related to the production efficiency of the manufacturing system. Research on these elements has received considerable attention from scholars. However, in most existing literature, these three elements are often studied individually or in pairs. Traditional research mainly focuses on the analysis of individual elements in manufacturing systems. For example, \cite{zhao2020two} proposed a two-stage cooperative evolution algorithm to solve a production scheduling problem in a non-waiting flow shop with the single consideration of production scheduling. \cite{cheng2020data} explored the way in which the internet promotes predictive maintenance, improving the feasibility of the assumed-based machine reliability maintenance strategy in manufacturing systems. Due to the close interaction among the elements in the manufacturing system, independent decision-making may lead to deviations from actual production situations. Following the demands, numerous scholars have attempted to collaboratively investigate the joint optimization of elements in the manufacturing system.

\subsection{\emph {Problem Formulation within a pair of elements}}
Current research in this field can be broadly categorized into two areas. The first category analyzes both production scheduling and reliability. Notable studies include \cite{aggoune2004minimizing}, who explored scheduling complexities in a process flow shop considering machine reliability, and \cite{uit2021joint}, who investigated machine reliability control through production rate adjustments. \cite{basciftci2020data} proposed a data-driven degradation model that incorporates maintenance and production scheduling while accounting for load dependence. It should be emphasized that in this model, production scheduling is regarded as independent of the machine's current operational status. A more recent line of research explores the interrelationship between a pair of the factors. \cite{kim2022optimal} conducted a study on the interrelationship between machine degradation rate and job processing time in a single-machine scheduling problem. \cite{briskorn2024scheduling} studied job-dependent and stochastic machine deterioration affected by processing jobs under a fixed job sequence.

The second category of research involves an analysis of both the quality as well as either job scheduling or reliability within a production system. \cite{chen2005quality} established the foundational coupling relationship between quality and reliability. \cite{lu2017opportunistic} proposed a novel maintenance strategy for multi-stage manufacturing systems that considered the interaction between machine failure rate and quality. Furthermore, various studies have analyzed the relationship between quality and production scheduling. \cite{eben2022relationships} developed a model that establishes the relationship between production scheduling and quality, proposing an integrated production and quality inspection system. \cite{azimpoor2021joint} proposed an integrated approach for production, inspection, and quality planning, taking into account the delay time failure process. These studies primarily focus on monitoring and proposing robust planning strategies to address uncertain non-conforming products. For instance, \cite{hao2015soft} devised a soft-decision-based approach to handle the uncertainty caused by quality considerations. \cite{sabri2022novel} developed a method to monitor the quality error rate and its impact on product reliability. In sophisticated manufacturing systems, the attainment of high-quality production necessitates simultaneous consideration of job scheduling, machine reliability maintenance, and product quality checks throughout the entire production process. Due to their interdependence and interaction, it is imperative to analyze these three key elements concurrently.

\subsection{\emph{Refining the problem formulation through triadic elements}}
In response to the escalating demand for production systems, recent research has incorporated a synthesis of three foundational elements, including quality, production scheduling, and machine reliability. For instance, \cite{cheng2020joint}  conducted a comprehensive study encompassing these three elements within a manufacturing system but simplified the job scheduling aspect by setting an index based on machine production efficiency and assuming a predefined deterioration path for the machine. \cite{tasias2022integrated} proposed a Bayesian model that ingeniously merged quality and machine reliability in a single-machine system, where the production process was presumed to exhibit a constant production efficiency. Similarly, \cite{rivera2021joint} treated production scheduling as a control variable in a manufacturing system, and \cite{shi2024optimising} assumed that products are produced at a constant rate in a imperfect manufacturing system. \cite{fakher2018integrating} integrated these three pivotal elements into a comprehensive model, wherein production scheduling predominantly contributes to the overall product quality, while disregarding the impact of job scheduling on production decisions within a manufacturing system. \cite{sinisterra2023delay} developed a model that incorporated the scheduling of a sequence of resumable jobs and inspection policies in a critical single-component system, where the occurrence of defects was contingent upon the hypothesized deterioration process. Reflecting on the aforementioned literature, it is evident that research encompassing quality, job scheduling, production, and reliability tends to simplify the production scheduling aspect by representing it solely as the machine production rate. Consequently, there is currently a lack of a more comprehensive analysis for coordinating job scheduling, quality inspection, and machine maintenance within the manufacturing system.

\subsection{\emph{Studies on uncertainties in rework-based production systems}}
In order to achieve a high-quality output in terms of both quantity and quality, various research efforts have incorporated rework activities into the production process \citep{sonntag2018disposal, berling2022inventory}. These studies have considered the existence of random defect rates within the production process. For instance, \cite{gouiaa2018integrated} developed an analytical model for maintenance decision-making, where non-conforming products undergo rework and the rate of quality deterioration is pre-assumed. In a similar vein, \cite{gautam2021integrated} assumed a uniform distribution for the proportion of non-conforming products and established a green supply chain model that integrates product recovery management. Some other researchers have examined the impact of machine deterioration on products, where defective items are either added to the regular production queue at the end of the processing time or reworked promptly \citep{nobil2020optimal}. 

However, uncertainties arise in real production systems where rework is influenced by the QRP-effect. These uncertainties relate to the status of pending jobs and the degradation path of machines, as well as the initial product quality. Proactive measures are necessary to address uncertainties in manufacturing systems. Researchers have developed algorithms employing fuzzy methods to proactively mitigate uncertainties in advance \citep{gao2020solving}. In addition, robust measures such as surrogate measures \citep{yang2020robust} and slack-based measures \citep{hazir2010robust} are commonly used to plan for uncertain production environments. Recently, flexible strategies have been introduced to dynamically schedule production and enhance efficiency. These flexible strategies utilize conditional information to adapt plans to uncertain production environments, outperforming static algorithms \citep{shi2020condition}. For example, a two-stage framework comprising offline and online stages, which incorporate \emph{in-situ} information, is frequently utilized in production problems. Online information including the work in progress and the buffer levels is employed to update the initial solution and improve its performance \citep{hoffman2021online}. Conditional job information is utilized to develop rescheduling strategies for updating the job arrangement \citep{wang2019improved}. However, a limitation of these solution frameworks is the lack of communication between the two stages, which hampers the enhancement of the initial solution's quality and restricts the performance of the updated solution by not leveraging valuable information from the digital platform.

In conclusion, when dealing with production systems that involve rework influenced by the QRP co-effect, uncertainties permeate the entire production process. Therefore, it is crucial to develop an adaptive approach that jointly determines maintenance strategies, production scheduling, and rework activities in order to effectively navigate and optimize the complex dynamics of such production systems.

\section{Problem Definition}

\subsection{\emph{Parallel Systems with Reworking Processes}}
We consider an unrelated parallel system with rework consideration under the QRP-co-effect, where non-conforming products are reprocessed to meet quantity demands. The general problem can be formally described as follows:
There are $n$ independent jobs, denoted as $N=\left \{ 1,2,...,n \right \}$, with different qualities, which need to be processed on $m$ distinct parallel machines, denoted as $M=\left \{ 1,2,...,m \right \}$. The manufacturing system aims to output $N$ conforming products. Due to machine degradation and variation in input job qualities, some non-conforming products require rework during the production process. These non-conforming products and their corresponding jobs are added to the list of pending jobs for reprocessing. To mitigate production losses caused by machine degradation, appropriate maintenance strategies, including corrective maintenance (CM) with replacement capability and imperfect preventive maintenance (PM), are implemented to maintain system reliability. The entire manufacturing system exhibits uncertainty and dynamic characteristics due to factors such as machine degradation paths, job processing times, pending job sequences (including rework jobs). Consequently, the joint optimization problem involves determining the optimal job scheduling, maintenance strategy, and quality control concurrently.


The joint optimization problem is based on the following assumptions:
\begin{enumerate}
  \item Initially, both jobs and machines are assumed to be in an available state.
  \item CM activity involves restoring a machine to its original operational status when it exceeds the failure threshold.
  \item Maintenance activity is not allowed while a machine is processing a job.
\end{enumerate}
\subsection{\emph{Modeling QRP-co-effect}}
Based on the analysis of the interrelationship between quality $Q_{i,k}$, machine reliability $R_{k}$, and scheduling $P_{i}$ in the production process, three intermediate variables with uncertainties need to be introduced: machine degradation, actual job processing time, and product quality. Considering the influence of the working environment on machine degradation state $W_{i,k}$ and the QRP-co-effect, the degradation path of the machine can be represented by:
\begin{flalign}
\Delta W_{i,k} = \Delta U_{i,k}^{(-)} + \Delta U_{i,k}^{(+)} + \Delta V_{i,k}, \label{6_3}
\end{flalign}
where $\Delta W_{i,k}$ denotes the change in machine degradation resulting from processing a job, while $U_{i,k}$ and $V_{i,k}$ represent the influence induced by the job and the environment, respectively. Specifically, $\Delta U_{i,k}^{(-)}$ represents the wear of machine $k$ after processing job $i$ whose initial quality is ineligible, following a Gaussian distribution \citep{ye2019reliability}: $\Delta U_{i,k}^{(-)}\sim Normal(\delta_{i}\mu_{k}^{(-)},\sigma_{k}^{(-)})$, where the product quality deviation $\delta_{i}=|D_{i}-SL_{i}^{(+)}|$, where $D_{i}$ and $SL_{i}^{(+)}$ represent the quality characteristic and the specification limit of the product corresponding to job $i$, respectively. $\Delta U_{i,k}^{(+)}$ represents the effect of job $i$ , and $\Delta U_{i,k}^{(+)}(p_{i}) \sim Normal(p_{i}\mu_{k}^{(+)},\sigma_{k}^{(+)})$, where $p_{i}$ represents the actual processing time of job $i$. $\Delta V_{i,k}$ illustrates the effect caused by the external environment, which is characterized by a Gamma distribution: $V_{i,k}(\Delta t_{i,k}) \sim Gamma(\alpha_{k}\Delta t_{i,k},\beta_{k})$, where $\Delta t_{i,k}$ denotes the time interval between the start time of job $i$ and the next job on machine $k$, expecting the maintenance time. 

The actual processing time of jobs is determined by the degradation state of the assigned machine and the nominal processing time, obtained as follows:
\begin{flalign}
p_{i,k}=O_{i,k}(1+\eta W_{i^{'},k}),\label{7_1}
\end{flalign}
where $p_{i,k}$ and $O_{i,k}$ represent the actual and nominal processing time of job $i$ processed on machine $k$, respectively. Job $i^{'}$ is the job processed prior to job $i$ on machine $k$. The actual processing time of the job increases with the deterioration of the assigned machine. The decline in machine condition leads to a decrease in the quality of processed jobs, resulting in an increased number of non-conforming products and the need for reworking activities. To ensure high sustainability and efficient operation of the system, an appropriate maintenance strategy is necessary to maintain machine’s condition. Production scheduling, which determines job sequencing, serves as the link between product quality and machine reliability. The impact of machine degradation on product quality is expressed as follows:
\begin{flalign}
\boldsymbol D=f(\boldsymbol \upsilon, \boldsymbol W),\label{8_3}
\end{flalign}
where $\boldsymbol D$ is a dimensional vector representing the quality characteristics of products under the degradation states of machines $\boldsymbol W$ and the initial job quality $\boldsymbol \upsilon$, with the characteristics of machines being independent of each other, as well as the initial jobs. In this production system, the output of non-conforming products from the manufacturing system can lead to delivery delays and necessitate rework, which also affects the assigned machine’ condition.
\section{Joint formulation with QRP-co-effect}

The decision variables in our formulation are summarized in Table \ref{decision}. The binary decision variables encompass job scheduling decisions, which involve determining the sequence of jobs and assigning them to specific machines. Additionally, maintenance decisions are incorporated, including the type of maintenance and the timing of its execution. Several continuous decision variables are defined, such as the start time of jobs and the reworking criteria, which are based on the current rate of non-conformance. Formulas (\ref{C1}-\ref{C5}) in the introduction section outline the framework of the model, while the subsequent subsections will provide detailed explanations and attributes of this integrated model.
\spacingset{0.9}
\begin{table}[t]
\centering
\caption{Decision variables for Modeling.}
\label{decision}
\renewcommand{\arraystretch}{0.9} 
\begin{tabularx}{\textwidth}{p{0.2\textwidth}lX}
   \toprule
       \emph{Job scheduling} & $x_{i^{'},i}\in \left \{0,1\right \}$ & 1, job $i$ is processed directly after job $i^{'}$ on the  same machine, and 0 otherwise; \\ 
         & $x_{i,k}\in \left \{0,1\right \}$ & 1, job $i$ is assigned on machine $k$, and 0 otherwise; \\
         & $S_{i}\geq 0$ & Start time of job $i$; \\ 
         & $C_{i}\geq 0$ & Completion time of job $i$. \\ 
       \midrule  
        \emph{Maintenance} & $z_{i}\in \left \{0,1\right \}$ & 1, the CM activity is performed after processing job $i$, and 0 otherwise; \\ 
         & $y_{i,g}\in \left \{0,1\right \}$ & 1, the PM activity belonging to the jointly maintenance group $g$ is performed after processing job $i$, and 0 otherwise; \\
         & $T_{g}^{pm}\geq 0$ & The maintenance time interval of group $g$ for the PM.\\ 
         \midrule  
        \emph{Quality} & $Thr^{r}\in[0,1]$ & The rescheduling decision point. \\
    \bottomrule
\end{tabularx}
\end{table}
\spacingset{1.5}
\spacingset{0.9}
\begin{table}[t]
\centering
\caption{Miscellaneous variables for problem formulation.}
\label{tab1}
\renewcommand{\arraystretch}{1.0} 
\begin{tabularx}{\textwidth}{p{0.1\textwidth} X}
\hline
\textbf{Symbol} & \textbf{Description} \\
\hline
$p_{i}$ & The actual processing time of job $i$; \\ 
$\Delta W_{i,k}$ & Deterioration amount under job-sequence affect during job $i$ processed on machine $k$; \\ 
$\Delta U_{i,k}^{(-)}$ & The non-conforming workload-effect deterioration of machine $k$ caused by job $i$; \\ 
$\Delta U_{i,k}^{(+)}$ & The qualified workload-effect deterioration of machine $k$ caused by job $i$; \\ 
$\Delta V_{i,k}$ & Deterioration amount under environmental noise during job $i$ processed on machine $k$; \\ 
$D_{i}$ & The quality characteristic of the product corresponding to job $i$;\\ 
$q_{i}$ & The record variable indicating whether the product corresponding to job $i$ is qualified;\\
$C_{max}$ & The makespan of the series-parallel multistage production system; \\ 
$C^{m}$ & The total cost of maintenance; \\ 
$C_{g}^{pm}$ & The maintenance cost of group $g$ for the PM. \\ 
\hline
\end{tabularx}
\end{table}
\spacingset{1.5}

\spacingset{1.0}
\begin{table}[t]
\begin{center}
\caption{Parameter setups in the formulation.}
\label{tab2}
\renewcommand{\arraystretch}{0.7} 
\begin{tabularx}{\textwidth}{lX}
\toprule
$N=\left \{ 1,2,...,n \right \}$, & The set of job indexes; \\ $K=\left \{ 1,2,...,k \right \}$, &  The set of machine indexes;\\ 
$N_{k}=\left \{ 1,2,...,n_{k} \right \}$ & The set of job indexes processed on machine $k$;  \\ 
$O_{i,k}$  & The nominal processing time of job $i$ scheduled on machine $k$;  \\ 
$\eta $ & The deterioration factor for the job actual processing time;  \\ 
$(T_{k}^{pm}, T_{k}^{ps}, C_{k}^{pm})$ & The preventive maintenance (PM) time duration, the associated setup time and PM cost for machine $k$; \\ 
$(T_{k}^{cm}, C_{k}^{cm})$ & The corrective maintenance (CM) time duration and the CM cost for machine $k$;   \\ 
$L_{k}$ & The CM threshold for machine $k$;\\ 
$(\alpha _{k}$, $\beta _{k})$ & The shape and the scale parameters of Gamma distribution for machine $k$ in degradation model;   \\ 
$(\mu_{k}^{(-)}$, $\sigma_{k}^{(-)})$ & The mean and the variance parameters of Gaussian distribution for machine $k$ in degradation model, where the processed job is non-conforming;   \\ 
$(\mu_{k}^{(+)}$, $\sigma_{k}^{(+)})$ & The mean and the variance parameters of Gaussian distribution for machine $k$ in degradation model, where the processed job is qualified; \\ 
$(\theta, \varphi)$ & The machine status influence factor and PM times influence factor in the imperfect PM policy;\\ 
$(\upsilon_{i}, \varepsilon_{i})$ & The initial quality of the job, and the vector of noise-variables for the quality characteristic;\\ 
$(\boldsymbol a$, $\boldsymbol b, \Gamma)$ & The coefficient vectors for the quality characteristics;\\
$(SL_{i}^{(+)}, \xi_{i})$ & The specification and the qualified threshold for the product corresponding to job $i$;\\
$(\mu, \sigma)$ & The mean and the variance parameters of Gaussian distribution for the characteristic of the job's initial  quality.\\
\bottomrule
\end{tabularx}
\end{center}
\end{table}
\spacingset{1.5}

\subsection{\emph{Modeling the QRP-co-effect}}
The QRP-co-effect of this problem can be addressed by formulating constraints within three main parts: job scheduling, maintenance, and quality considerations. In the first part, decisions concerning the assignment of machines to jobs and the sequencing of jobs on the machines need to be made. Furthermore, the impact of machine reliability on production scheduling is also considered. The second part encompasses constraints related to maintenance activities. This includes characterizing the constraints among different maintenance actions. Additionally, the machine's degradation path, changes in product quality, and the actual processing time of jobs under the QRP-co-effect during production are formulated. Lastly, in the third part, constraints regarding the relationship between the quality states of products and the rescheduling point of the production system are provided.

\subsubsection{\emph{Constraints for job scheduling}}
Within the production system, each job is assigned to a single machine and can only be processed on that machine. Additionally, for any given machine, only one job (except for the last job) can be processed consecutively after another job, denoted as job $i^{'}$. The specific conditions for these constraints are expressed as follows:
\begin{alignat}{2}
\label{C4} & \sum_{k\in K}x_{i,k}=1, \quad \forall i\in N_{k},&\\
\label{C5} & x_{i^{'},i}+x_{i,i^{'}}\leq 1, \quad \forall i\in N_{k},&\\
\label{C6} & \sum_{i^{'}\in N_{k}}x_{i^{'},i}\leq 1,\quad \forall i\in N_{k},&\\
\label{C7} & \sum_{i\in N_{k}}x_{i,k}x_{0,i}= 1, \quad \forall k \in K.&
\end{alignat}
Constrain (\ref{C7}) enforces that the dummy job 0 can only serve as a predecessor to one job on each machine. Additionally, inequality (\ref{C8}) restricts the start time of job $i$ to be greater than or equal to the completion time of the previous job $i^{'}$ on the same machine, considering the cost of maintenance time. Moreover, the completion time of a job must be greater than the sum of its start time and the actual processing time, as governed by condition (\ref{C9}).
\begin{alignat}{2}
\label{C8} 
S_{i} - S_{i^{'}}  &\ge p_{i^{'}}+\sum_{g\in G}y_{i^{'},g}T_{g}^{pm} +\sum_{k\in K}x_{i^{'},k}\sum_{i^{'}\in N_{k}}z_{i^{'}}T_{k}^{cm}  \nonumber \\
&\quad -(1-x_{i^{'},i})M, \quad \forall i\in N, i^{'}\in N, &\\
\label{C9} & C_{i} \geq  S_{i}+ p_{i},  \quad \forall {i} \in {N}, &\\
\label{C10} & x_{i^{'},i}\in \{0,1\}, \quad \forall i\in N,i^{'}\in N,i\neq i^{'} &\\
\label{C11} & x_{i,k}\in \{0,1\}, \quad \forall i\in N, k\in K,i\neq i^{'} &\\
\label{C12} & S_{i}>0, \quad \forall i\in N.&
\end{alignat} 
Constraints (\ref{C10}-\ref{C12}) determine the range of the decision variables for the job assignment.

\subsubsection{\emph{Constraints for machine maintenance}} 
Constraint (\ref{C13}) guarantees that consecutive or simultaneous maintenance activities are not allowed. Constraint (\ref{C14}) addresses the concept of maintenance grouping, where $T_{g}^{pm}$ represents the duration of the preventive maintenance activity, and $y_{i,g}$ is an intermediate variable defined in (\ref{C18}) to determine whether maintenance grouping should be executed. Constraints (\ref{C15}-\ref{C16}) ensure that each preventive maintenance activity is assigned to a single maintenance group, and prevent two preventive maintenance events for a machine from belonging to the same maintenance group. Finally, constraint (\ref{C17}) specifies that corrective maintenance events must be promptly carried out upon machine failure. Constraint (\ref{C18}) indicates whether the PM should be executed.
\begin{alignat}{2}
\label{C13} & \sum_{g\in G}y_{i,g}+z_{i}\leq 1, & \forall i\in N,\\
\label{C14} & T_{g}^{pm}\geq \sum_{i\in N_{k}}x_{i,k}y_{i,g}(T_{k}^{pm}+T_{k}^{ps}),& \forall k\in K,  \\
\label{C15} & \sum_{g\in G}y_{i,g}\leq 1, & \forall i \in N,\\
\label{C16} & \sum_{i\in N_{k}}x_{i,k}y_{i,g}\leq 1,& \forall  k\in K, g\in G,\\
\label{C17} & \begin{array}{l}
z_{i}=\left\{\begin{matrix}0, \ \sum_{k\in K}x_{i,k}(W_{i,k}-L_{k})<0,\\ 
1, \ \sum_{k\in K}x_{i,k}(W_{i,k}- L_{k})> 0,
\end{matrix}\right. 
\end{array} &  \forall i\in N,\\
\label{C18} & y_{i,g}\in \{0,1\},  &  \forall i\in N, g\in G. &
\end{alignat}

\subsubsection{\emph{Characterization for machine deterioration, quality and production scheduling dynamics in production systems}}

(1) \emph{The deterioration of each machine}

The deterioration state of the machine depends on the job scheduling and the maintenance activity, and the external environmental factors, which can be expressed as follows:
\begin{alignat}{2}
& 
W_{i,k}=(1-\sum_{i^{'}\in N_{k}}x_{i^{'},i}z_{i^{'}})[I(\sum_{i^{'}\in N_{k}}x_{i^{'},i}  \nonumber \\
&\quad (W_{i^{'},k}+\Delta W_{i,k})| \sum_{i^{'}\in N_{k}}x_{i^{'},i}y_{i^{'},g})], &  \forall i\in N_{k}, k \in K, & \label{W1}
\end{alignat}
where function $I(d|e)$ can be represented as follows:
\begin{alignat}{2}
\label{C23} & I(d|e)=\left\{\begin{matrix}
d,&e=0  \\ 
\theta  d+\varphi N_{i,k}^{pm},&e=1 
\end{matrix}\right.
.&
\end{alignat}
In formula (\ref{C23}), term $\theta x + \varphi N_{i,k}^{pm}$ describes the impact of imperfect preventive maintenance (PM) on the machine's state $W_{i',k}$ and the number of PMs performed since the most recent corrective maintenance (CM) or the startup time of the machine $N_{i,k}^{pm}$. Details for obtaining $N_{i,k}^{pm}$ can be seen in Appendix A. The change $\Delta W_{i,k}$ for machine $k$ under job-sequence affect during job $i$ in equation (\ref{W1}) can be calculated as follows:
\begin{alignat}{2}
\label{C22} & 
\Delta W_{i,k}= \Delta U_{i,k}^{(+)}(p_{i})+\Delta U_{i,k}^{(-)}(\delta_{i})+\Delta V_{i,k}(\Delta t_{i}), \forall i\in N_{k}, k \in K, & 
\end{alignat}
where intermediate variable $\Delta t_{i}$ can be obtained by:
\begin{alignat}{2}
\label{C21} & \Delta t_{i}=S_{i}-\sum_{i^{'}\in N}x_{i^{'},i}S_{i^{'}} -\sum_{i^{'}\in N}x_{i^{'},i} \nonumber \\
&\quad (\sum_{g\in G}y_{i^{'},g}T_{g}^{pm}+ \sum_{k\in K}x_{i^{'},k}\sum_{i^{'}\in N_{k}}z_{i^{'}}T_{k}^{cm}), & \forall i\in N.&
\end{alignat}
In formula (\ref{C21}), $\Delta t_{i}$ represents the time interval (excluding maintenance time) between the start time of the production of job $i'$ and the start time of the production of job $i$. 

(2) \emph{Product quality}

The quality of the output product is influenced by both the degradation state of the machine assigned through production scheduling and the initial quality of the input job. The quality attribute of product $i$ can be calculated using the following equation:
\begin{flalign}
D_{i}=\sum_{k\in K}x_{i,k}(\upsilon_{k}^{0} + a_{k} \sum_{i^{'}\in N}x_{i^{'},i} W_{i^{'},k}+ b_{k}^{0} \cdot  \varepsilon_{i}+\sum_{i^{'}\in N}x_{i^{'},i}W_{i^{'},k} \Gamma_{k} \cdot  \varepsilon_{i}),\label{C29}
\end{flalign}
where $\upsilon_{i}=x_{i,k}\upsilon_{k}^{0}$ represents the initial quality of the job. Parameters $b_{i}=x_{i,k}b_{k}^{0}$, $a_{k}$, $v_{k}$ and $\Gamma_{k}$ are impact vectors. $\varepsilon_{i}$ denotes the vector of noise-variables. Here the characteristic of the job's initial quality is assumed as following a Truncated normal distribution $N(\mu, \sigma)$.

where $\boldsymbol D$ is a dimensional vector representing the quality characteristics of products under the degradation states of machines $\boldsymbol W$, with the characteristics being independent of each other. $\boldsymbol \upsilon$ is the initial job quality vector following a truncated normal distribution. $\boldsymbol \varepsilon$ denotes the vector of noise variables. $\boldsymbol a$ and $\boldsymbol b$ are impact vectors that describe the effect of machine degradation states and noise on the products, respectively. $\boldsymbol \Gamma$ characterizes the interaction effect between machine degradation $\boldsymbol W$ and noise $\boldsymbol \varepsilon$. In this production system, the output of non-conforming products from the manufacturing system can lead to delivery delays and necessitate rework, which also affects the assigned machine’ condition.

The quality of the product is assessed based on the specification $SL_{i}^{(+)}$ and the qualified threshold $\xi_{i}$ for the characteristic of job $i$:
\begin{alignat}{2}
\label{C30}&q_{i}= \left\{\begin{matrix}
1, \quad |D_{i}-SL_{i}^{(+)}|<\xi_{i}\\ 
0, \quad |D_{i}-SL_{i}^{(+)}|\geq \xi_{i},
\end{matrix}\right.
\end{alignat}
where $q_{i}$ is used to record the qualified products.


(3) \emph{Job actual processing time}

The actual processing time of the job in production is influenced by the deteriorating condition of the assigned machine, as indicated by formula (\ref{7_1}):
\begin{flalign}
\label{C34} p_{i} = \sum_{k\in K}x_{i,k}O_{i,k}(1+\eta W_{i,k}),  \forall i\in N. &
\end{flalign}

(4) \emph{Reworking for non-eligible products}

We employ a partial rescheduling strategy, which is triggered when a non-conforming rate reaches a certain threshold:
\begin{flalign}
\sum_{i \in N_{l}}(1-q_{i})/|N_{l}| \geq Thr^{r} ,&\label{ru_1}
\end{flalign}
where $N_{l}$ represents the set of jobs processed from the start of scheduling or since the $(l-1)^{th}$ rescheduling point.

A fitness function is formulated as follows to evaluate the rescheduling scheme, considering both the time cost and maintenance cost associated with processing unit jobs: 
\begin{flalign}
f^{r}=((\sum_{i \in N_{l}}q_{i})^{2}/(C_{l}^{m,r}T_{l}^{r}), \quad \forall l\in L.&
\end{flalign}
Here, $\sum_{i \in N_{l}}q_{i}$, $C_{l}^{m,r}$ and $T_{l}^{r}$ represent the number of conforming products, maintenance cost, and time cost generated during the $l^{th}$ rescheduling in the production process, respectively. 

\subsection{\emph{Objective function}}
The objective of this study is to maximize overall productivity while minimizing maintenance costs. This is determined by two key factors: the makespan $C_{max}$ of jobs and the overall maintenance cost $C^{m}$. The objective function can be expressed as follows:
\begin{alignat}{2}
\label{C36} & C_{max}=max (C_{i}),  & \ \ \forall i\in N,\\
\label{C37} & C^{m}=\sum_{i\in N}z_{i}\sum_{k\in K}x_{i,k}C_{k}^{cm}+\sum_{g\in G}C_{g}^{pm},\\
\label{C38} & \begin{array}{l}
C_{g}^{pm}=\sum_{i\in N_{k}}x_{i,k}y_{i,g}(C_{k}^{pm}+C_{k}^{ps}), \quad
\end{array}& \quad \forall k\in K, g\in G,
\end{alignat}
where the setup cost of preventive maintenance (PM) is distributed evenly among the machines within the joint maintenance group $g$.

\section{Dual-module algorithm}
This integrated optimization problem considering rework, is known to have NP-hard computational complexity \citep{hu2022scheduling}, thereby posing a substantial computational challenge, particularly when the number of jobs and machines involved grows significantly. Deterministic optimization methods may encounter difficulties in efficiently optimizing the objective function, both in terms of time and due to inherent uncertainties.

However, leveraging the capabilities of a data platform enables the acquisition of real-time information on machine status and rework progress after job task execution. This, in turn, allows for dynamic updates to the existing plan. In this section, an algorithmic framework called \underline{D}ual-module \underline{P}lanning and \underline{E}valuation \underline{I}ntegration \underline{A}lgorithm (DPEIA) is introduced to address this integrated optimization challenge. The goal of DPEIA is to facilitate more effective and integrated decision-making in scheduling and maintenance, ultimately enhancing overall efficiency. Compared to the traditional two-stage approach where planning and scheduling are performed independently, DPEIA introduces a dynamic dual-module communication mechanism to enhance the planning efficiency by actively interacting with the simulation environment. This mechanism also improves computational efficiency by reducing the need for frequent rescheduling caused by non-conforming products.

\begin{figure}[t]
\centering
\includegraphics[scale=0.6]{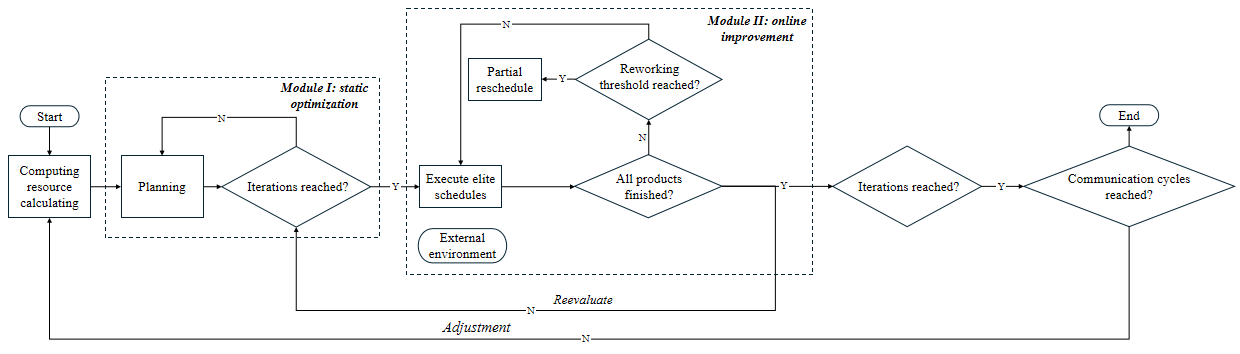} \label{algorithm}
\caption{\textcolor{black}{The framework of the dual-module algorithm.}}
\label{Figure_algrithmframe}
\end{figure}

The dual-module algorithmic framework, as illustrated in Figure \ref{Figure_algrithmframe}, consists of a planning module (static optimization without rescheduling) and an evaluation module (online improvement with rescheduling). The planning module aims to obtain a stable production scheduling plan with computational efficiency, while the evaluation module dynamically adjusts the search direction of the planning module through communications between the modules. Specifically, elite individuals are generated through static optimization and undergo practical reevaluation by interacting with the evaluation module. This guides the search direction of the planning module in the next iteration. Multiple rounds of static optimization and online improvement are performed in an interactive session to obtain an ultimate plan. Furthermore, this section proposes an adaptive maintenance strategy based on system condition information, such as machines' degradation status, job pending list, and product quality, among others.

The following subsections will provide detailed explanations of the DPEIA method.

\subsection{\emph{Static optimization}}
In the planning module, our objective is to generate a resilient plan for job scheduling that includes idle spaces, without the need for rescheduling. This framework is adaptable and can accommodate various optimization algorithms. We have chosen to incorporate an \underline{E}nhanced version of the \underline{M}ulti-\underline{O}bjective \underline{D}ifferential \underline{E}volution algorithm (EMODE), due to its flexible search capability. To address the impact of uncertain rework activities on the schedule, idle spaces are introduced as a mitigation strategy. The fitness function used in the static optimization module can be defined as follows:
\begin{alignat}{2}
f^{s}=((\sum_{i \in N}q_{i})^{2}/(C^{m}C_{max}).\label{fit_1}
\end{alignat}
This evaluation is conducted by minimizing both the time cost and maintenance cost related to processing unit jobs.

In this module, the quantity of idle space $n_{e}$ is obtained using formula (\ref{n_1}-\ref{n_2}), given that the impact of idle space on machine degradation is equivalent to that of qualified jobs on the machine's degradation path. To accommodate the differing processing capabilities of machines for various types of jobs, the jobs are divided based on their types. Specifically, jobs of the same type are assigned to machines with corresponding processing capabilities. All machines are initially set to their initial states, denoted as $\boldsymbol W^{0}$. The total number of non-conforming products $n^{0}_{e}=\sum_{i \in N}(1-q_{i})$, considering the machine's degradation state $W_{i,k}^{0}$. 

The number of idle spaces $n_{e^{th}}$ caused by jobs of the $e^{th}$ type can be expressed as:
\begin{flalign}
n_{e^{th}}=n^{0}_{e^{th}}/M_{e^{th}},
\label{n_1}\\
n_{e}=\sum_{e^{th}}n_{e^{th}},\label{n_2}
\end{flalign}
where $M_{e^{th}}$ denotes the number of machines capable of processing jobs of the $e^{th}$ type. The efficiency of this algorithm can be further enhanced based on the following proposition.

\emph{Proposition 1}: For the proposed problem, If two consecutive jobs processed on the same machine meet any of the following conditions, their processing sequences should be interchanged:

(1) Both output products are non-conforming, and the nominal processing time of the first job $i^{'}$, exceeds that of the second job $i$.

(2) Both output products belong to the same type and are conforming. The nominal processing time of the first job $i^{'}$ exceeds that of the second job $i$ and satisfies the following condition:
\begin{flalign}
\Delta W_{i,k}(a_{k}+\Gamma_{k} \varepsilon_{i})<\Delta D_{i^{'}}.
\end{flalign}

Here, $\Delta W_{i,k}$ represents the difference of machine's states between two solution after finishing these two jobs. Furthermore, $\Delta D_{i^{'}}$ represents the discrepancy between the quality feature of job $i^{'}$ and the specification $SL_{i}^{(+)}$. The proof is represented in Appendix B.

In the production process, job processing plays a pivotal role in the machine's degradation process, encompassing the impact of factors such as job sequence, quality, and quantity on the machine's degradation trajectory. The EMODE algorithm incorporates two crucial differential evolution operators: the decoding-based Differential Evolution (DE) operator \citep{zhou2019self} and the similarity-oriented Ranking Evolution (RE) operator \citep{hou2023multi}. The DE operator introduces a significant amount of randomness, contributing to the enhancement of solution diversity. The RE operator takes advantage of the inherent data associations among different jobs, refining solutions to optimize the process. The control parameter $\nu$ in the EMODE algorithm is used to allocate computational resources for algorithm exploration and exploitation. In the early stages of the iterative process, the algorithm focuses on exploration to avoid getting trapped in local optima. As the process continues, the algorithm gradually shifts its emphasis to exploitation, leading to faster convergence. The algorithm performs well, particularly when the control parameter $\nu$ decreases over the iterations. This relationship is mathematically expressed as follows:
\begin{flalign}
\nu = 2(1-(\frac{\varrho }{\varrho^{0}})),\label{a_1}
\end{flalign}
where $\varrho^{0}$ represents the total number of iterations, and the current iteration number $\varrho \in[0, \varrho^{0}]$. In each iteration, the control parameter $\nu$ is randomly selected within the range $[0,2]$ using this function and $rand$ is utilized to generate a random number within the range (0,1).
In this algorithm, when $\nu>1$ and $rand<0.7$, the similarity-oriented RE operator is used to balance the number of pending jobs on each machine. 
 Otherwise, the DE operator is utilized to generate a new solution. To enhance the ability of local search, two jobs are randomly selected from an individual. If one of the conditions in \emph{Proposition 1} is satisfied, the processing sequence of the jobs will be swapped. Otherwise, a job will be randomly selected from the machine with the longest total work completion time and moved to the most idle machine. This algorithm is detailed in Algorithm \ref{Algorithm 1}(presented in Appendix G), where the idle spaces are randomly assigned to the machine at the initial stage. A planning solution can be obtained through the EMODE algorithm to provide a job allocation strategy.

\subsection{\emph{Online Improvement through Environmental Interaction}}
In the context of an online improvement module that incorporates rescheduling, each elite individual undergoes a reevaluation by interacting with the actual production environment. A partial rescheduling strategy is employed, which relies on conditional information to determine the rescheduling point. This point is determined based on the rate of non-conforming products, denoted as $Thr^{r}$ in formula (\ref{ru_1}). The determination of the decision rescheduling scope relies on the rescheduling strategy as presented in the referenced work \citep{wang2019improved}, involving the aggregation of available idle spaces on machines within the current planning module to accommodate the inclusion of a new job. In this module, the evaluation function is designed as follows:
\begin{flalign}
\label{C58}&f^{(eva)}=(1/(\sum_{l\in L}C_{l}^{m,o}\sum_{l\in L}T_{l}^{o}d^{o}),
\end{flalign}
where $\sum_{l\in L}C_{l}^{m,o}$ and $\sum_{l\in L}T_{l}^{o}$ represent the time and maintenance cost incurred for scheduling and maintenance in the production process, respectively, and $L$ represents the number of rounds re-scheduled during the actual production. Variable $d^{o}$ represents the deviation between the planned and actual scheduling in terms of job start times and assignments, with further details available in reference \citep{an2022multiobjective}. 

For each rescheduling, the promising region always exists in the neighborhood of the initial plan when new jobs are inserted in the unrelated parallel machines system \citep{peng2019multi}. Algorithm \ref{Algorithm 2} outlines a fast local search method, which utilizes the JS (Job Swapping) and JI (Job Insertion) operators to search for high-quality solutions \citep{ulaga2022local}. This method is particularly effective in scenarios where new jobs are inserted into an unrelated parallel machines system. The combination of JS and JI operators has been shown to yield favorable results. In this approach, a job is chosen with a 0.5 probability to be inserted at position $j$ on the machine with the earliest completion time. This condition is met when the nominal processing time of the selected job is shorter than the processing time of the job assigned at the $(j+1)^{th}$ position, but longer than the job assigned at the $(j-1)^{th}$ position. Alternatively, with a probability of 0.5, the positions of randomly selected jobs are swapped. We introduce \emph{Lemma 1} to guide the job scheduling decision.

\emph{Lemma 1}: In a parallel system, there always exists an optimal schedule where no idle time exists between two consecutive jobs on a machine.

The proof of \emph{Lemma 1} is provided in Appendix C.

Based on \emph{Lemma 1}, idle spaces in the scheduling are filled either by rescheduled jobs or pending jobs as per the plan. An example is given to illustrate this decoding process, which is represented in Figure \ref{Figure_decode}, provided in Appendix D. 
\subsection{\emph{Communication scheme and the enhanced maintenance decision}}
\subsubsection{\emph{Communication scheme}}
To address the challenge of obtaining a high-quality solution initially and the resource-intensive nature of high-frequency rescheduling, a communication scheme is designed to allocate computing resources between the static optimization module and the online improvement module. This scheme involves dividing the entire solving process into $N^{s}$ interactive rounds or communication cycles, with $maxIter$ representing the total computing resource. The allocation of computational resources for the two modules follows a Gaussian distribution pattern, as depicted in Figure \ref{percent} in Appendix F, which adjusts the computational proportions of the static planning module and the online improvement module in each communication cycle. The proportion of computational resources allocated to the online improvement module increases over time. In each round, the top $popSize^{P}$ individuals are selected based on the estimation of the function $f^{s}$ value, using $maxIter_{r}$ iterations from the previous cycle in the planning module. These selected individuals then undergo $maxIter^{R}_{r}$ iterations in the online improvement module to obtain a new label through function $f^{(eva)}$. This new label replaces the original $f^{s}$ in the static module and adjusts the search direction for the next round. The parameters $maxIter^{s}_{r}$ and $maxIter^{R}_{r}$ are updated at the end of each interaction, and the interaction is repeated $N^{s}$ times. Details of the computational resource allocation for each cycle are provided in Appendix F.

\subsubsection{\emph{Maintenance decision powered by DPEIA}}
To improve the effectiveness of maintenance decisions, we introduce a proposition within the framework of the DPEIA. Here, the term ``life cycle" of a machine refers to the time interval from its initial state to the occurrence of failure. To refine the maintenance decision, we then derive \emph{Proposition 2} of the above optimization problem as follows.

\emph{Proposition 2} : Considered both jobs and maintenance tasks, an optimal schedule exists where the PM activity should be suspended until machine failure occurs, given the following conditions:
\begin{flalign}
\frac{n^{pm}}{n^{c}}<\varsigma^{m},
\end{flalign}
where $\varsigma^{m}=min\{\frac{T^{pm}}{T^{c}}, \frac{C^{pm}}{C^{c}}\}$, and $n^{pm}$ represents the additional number of processed jobs on the machine when performing the PM task compared to not performing it during the current life cycle. $n^{c}$, $T^{c}$, and $C^{c}$ denote the number of processed jobs, job processing time, and maintenance cost, respectively, when the PM task is not executed on the machine in the current life cycle. The proof can be found in Appendix E.

The ``PM-suspension period" refers to the machine not performing the PM activity from the current point until it fails. In the DPEIA's maintenance decision process, a machine is classified as being in the ``PM-suspension period" if the condition stated in \emph{Proposition 2} is satisfied. No PM event should be carried out during the machine's current life cycle. Additionally, machines in this period will be screened out before joint maintenance is performed. Please refer to Appendix K for a detailed explanation of the algorithm's encoding and decoding components.
\section{Numerical Studies}
\subsection{\emph{Experiment setups}}
This section presents numerical experiments to showcase the effectiveness of the proposed approach, encompassing the maintenance strategy and dynamic exchange module within the production process. A heterogeneous parallel production system is considered, with careful consideration given to rework in order to address potential defects and maintain exceptional quality. The experiment focuses on the utilization of four unrelated parallel processing machines with distinct capacities, taking into account different initial statuses $\boldsymbol W^{0}$. Additionally, two types of jobs are considered, characterized by different nominal processing times. We follow \cite{ye2019reliability} to initiate the production conditions as follows:
The first type of jobs can be processed on machines 1, 3, and 4, with corresponding nominal processing times following a uniform distribution: $U\left(2.616, 0.3\right)$. The second type of jobs can be processed on machines 2 and 4, with nominal processing times following a uniform distribution: $U\left(1.92, 0.5\right)$. The quality of input jobs is assumed to follow a truncated normal distribution: $N(\mu_{q},\sigma_{q})$. The experiment explores three instances of job sizes, consisting of 100, 200, and 300 jobs, respectively. Table \ref{tablepa_3}, presented in Appendix I, outlines the parameter settings for this problem, where $SL^{(+)}_{e^{th}}$ and $\xi_{e^{th}}$ denote the specification limit and qualified threshold for the quality characteristic of the $e^{th}$ job.

In order to assess the effectiveness of the proposed method, four comparative algorithms are utilized: EMA (Enhanced Memetic Algorithm)  \citep{afsar2022multi}, MOEA-D (Multiobjective Evolutionary Algorithm based on Decomposition) \citep{de2022diversity}, MVO (Multi-Verse Optimizer) \citep{kumar2023two}, and NSGAII (Non-dominated Sorting Genetic Algorithm-II) \citep{lopez2021non}. The maintenance strategy employed by these comparative algorithms is based on the MCIM approach \citep{chen2022maintenance}, where the maximum number of preventive maintenance actions for new machines is obtained through static optimization, and a right-shift rescheduling policy is adopted. In evaluating the relative performance of the algorithms, three metrics are employed. The IGD-metric (Inverted Generational Distance) and hypervolume-metric are used to assess the approximation and distribution quality of the non-dominated solution set. Additionally, the RPD-metric (Relative Percentage Deviation) is utilized to analyze the maximum completion time $C_{max}$ and maintenance cost aspects of the results. Details of the metric calculations and algorithm parameter settings are provided in Appendix I.

\subsection{\emph{Performance comparison with benchmarks}}
The results of the algorithm comparison experiment based on the IGD-metric and hypervolume-metric are presented in Table \ref{table igd-hv}. Smaller values of the IGD-metric and larger values of the hypervolume-metric indicate more competitive algorithms. It is evident from the results that the DPEIA algorithm consistently achieves smaller objective function values across all instances compared to other algorithms, indicating its superior performance. Additionally, the DPEIA algorithm outperforms other algorithms with the maximum hypervolume value. Table \ref{tablerpdres} shows the experimental results of the DPEIA algorithm and the four other algorithms in terms of two objectives: maximum completion time and total maintenance cost, under the RPD-metric. Smaller values indicate better algorithm performance. The DPEIA algorithm demonstrates the best performance in both dimensions of maximum completion time and total maintenance cost.
\spacingset{0.9}
\begin{table}[t]
	\caption{{Comparison of algorithms performance based on various evaluation criteria.}}
	\label{table igd-hv}
         
	\setlength{\tabcolsep}{1.81mm}{
		\begin{tabular}{cccccccccc}
			\toprule
			Assessment criteria & Job sizes & $N$   & EMA & MOEA-D & MVO & NSGAII & DPEIA \\ 
			
			\midrule
			\multirow{3}{*}{IGD} & small-sized & 100 & 0.74 & 0.41 & 0.54 & 0.54 & 0.04\\
			& medium-sized & 200 & 1.07 & 0.85 & 0.93 & 0.99 & 0.03  \\
			& large-sized & 300 & 1.18 & 0.91 & 1.01 & 1.06 & 0.01  \\
			
			\midrule
			\multirow{3}{*}{HV} & small-sized & 100 & 0.11 & 0.23 & 0.18 & 0.18 & 0.49 \\
			& medium-sized & 200 & 0.13 & 0.19 & 0.13 & 0.10 & 0.49  \\
			& large-sized & 300 & 0.01 & 0.04 & 0.02 & 0.00 & 0.27  \\
			\bottomrule
		\end{tabular}
	}
\end{table}
\spacingset{1.5}

\spacingset{0.8}
\begin{table}[t]
	\caption{{Comparison of algorithms performance based on the RPD-metric.}}
    \label{tablerpdres}
    \setlength{\tabcolsep}{2.18mm}
    \begin{tabular}{ccccccccccc}
        \toprule
        \multirow{2}{*}{$N$} & \multicolumn{2}{c}{EMA} & \multicolumn{2}{c}{MOEA-D}  & \multicolumn{2}{c}{MVO}& \multicolumn{2}{c}{NSGAII}
        & \multicolumn{2}{c}{DPEIA}\\
        \cmidrule(r){2-3}\cmidrule(r){4-5}\cmidrule(r){6-7}  \cmidrule(r){8-9}  \cmidrule(r){10-11}  
         & $C_{max}$ & $C^{m}$ & $C_{max}$ & $C^{m}$ & $C_{max}$ & $C^{m}$ & $C_{max}$ & $C^{m}$ &
        $C_{max}$ & $C^{m}$\\
        \midrule
        100 & 98.18 & 258.71 & 88.48 & 156.45 & 55.37 & 198.04 & 63.32 & 211.83 & 4.43 & 45.72\\
        200 & 96.17 & 99.62 & 70.96 & 85.56 & 68.67 & 102.77 & 70.65 & 100.38 & 3.39 & 4.63\\
        300 & 121.93 & 73.10 & 92.01 & 69.03 & 81.45 & 75.23 & 67.39 & 83.06 & 3.24 & 1.84\\
        \bottomrule
    \end{tabular}
\end{table}
\spacingset{1.5}


\begin{figure}[t]
\begin{minipage}[b]{0.49\linewidth}
	\begin{center}
    	\includegraphics[width=\textwidth]{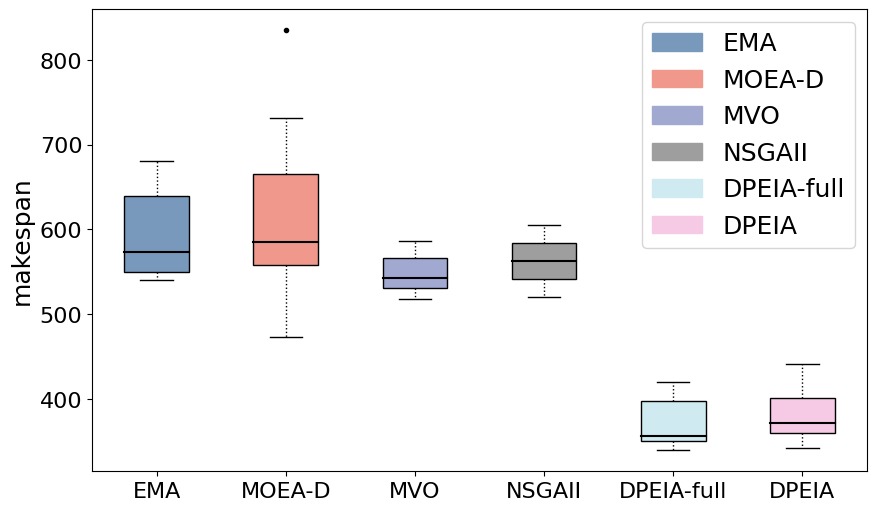}
     \caption{The box plot of makespan for 200 jobs under basic standard deviation.}
        \label{Figure_box_makespan}	
 \end{center}
\end{minipage}
  \hfill
\begin{minipage}[b]{0.49\linewidth}
	\begin{center}
    	\includegraphics[width=\textwidth]{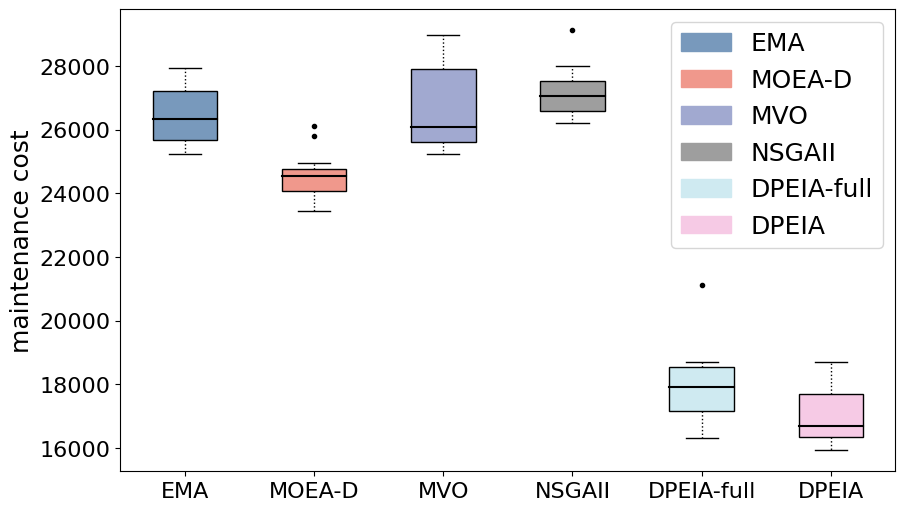}
     \caption{The box plot of maintenance cost for 200 jobs under basic standard deviation.}
        \label{Figure_box_totalcost}	
 \end{center}
\end{minipage}
\end{figure}

\begin{figure}[t]
\begin{center}
\begin{minipage}[b]{0.49\linewidth}
	
    	\includegraphics[width=\textwidth,
     height = 5.cm]{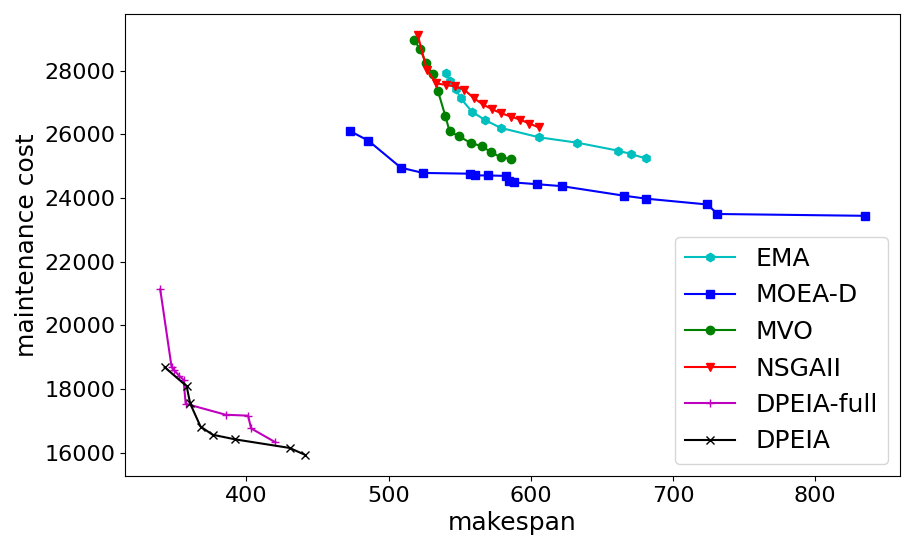}
     \caption{The Pareto optimal surface of maintenance cost for 200 jobs under basic standard deviation.}
        \label{Pareto_new}	
\end{minipage}
  \hfill
\begin{minipage}[b]{0.49\linewidth}

    	\includegraphics[width=\textwidth, height = 4.9cm]{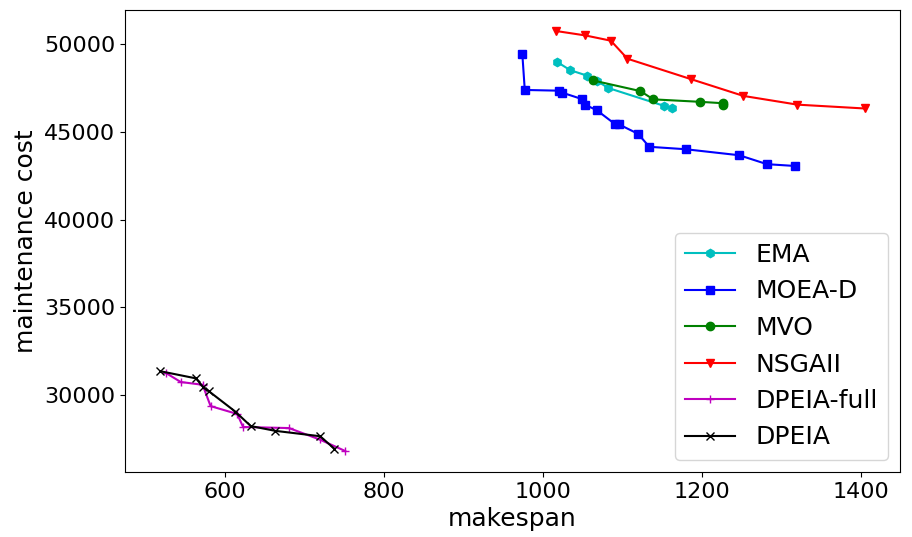}
     \caption{The Pareto optimal surface of maintenance cost for 300 jobs under the deviation of 0.09.}
        \label{Pareto_new1}	
 
\end{minipage}
\end{center}
\end{figure}

To provide a clearer visualization of the experimental results, box plots of the four other benchmarks, along with three algorithms of ablation experiments, are shown in Figure \ref{Figure_box_makespan} and Figure \ref{Figure_box_totalcost}, representing makespan and cost, respectively. The DPEIA and DPEIA-full methods show notable advantages in terms of makespan compared to the remaining four algorithms, highlighting the effectiveness of the proposed solution framework in uncertain manufacturing systems. Regarding maintenance cost, both DPEIA and DPEIA-full exhibit superior performance with minimal variance fluctuations compared to the other four algorithms, showcasing the advantages of simultaneous maintenance and production scheduling.

Figure \ref{Pareto_new} and Figure \ref{Pareto_new1} display the final set of non-dominated solutions obtained by DPEIA-full, DPEIA, and the four other algorithms. The solution set generated by the DPEIA method exhibits superior approximation and distribution diversity compared to the other algorithms. Furthermore, the DPEIA method achieves results that are very close to those obtained by the DPEIA-full algorithm while consuming approximately half of the computational resources. Based on the analysis of the experimental results, we conclude that DPEIA demonstrates exceptional performance in addressing the studied problems.

\subsection{\emph{Sensitive analysis}}
To validate the effectiveness of the proposed DPEIA algorithm, additional numerical experiments were conducted to analyze the sensitivity of the algorithm to different ranges of input job quality. Two qualified deviations, 0.03 and 0.09, were considered in the experiments. The results of these experiments, presented in Table \ref{tablerandom5}, can be observed that the DPEIA algorithm remains competitive compared to the other four algorithms in this experiment.

In addition, ablation experiments were conducted to evaluate the performance of the communication mechanism and the integrated mechanism separately. The DPEIA-full algorithm, where the proposed communication mechanism was replaced by a fully rescheduling strategy, and the DPEIA-single algorithm, which considers scheduling and maintenance decisions separately, were used as comparative benchmarks. Detailed comparison results between the DPEIA algorithm and the DPEIA-full algorithm in terms of IGD and hypervolume indicators are shown in Table \ref{tablerandom8}. 

The experimental results show that as the number of jobs or non-conforming input jobs increases, along with an increased frequency of maintenance demands, the disparities in cost and completion time among the algorithms become more pronounced. This emphasizes the advantages of the proposed method in dynamic production environments and its suitability for making maintenance decisions. Based on the IGD and hypervolume indicators, when the job qualified deviation $\sigma_{q}$ is 0.03, the performance of DPEIA is almost the same as that of DPEIA-full, while it is slightly worse when $\sigma_{q}=0.06$. DPEIA-full has advantages when $\sigma_{q}=0.09$, but situations with a high proportion of non-conforming jobs (about 45\%) are rare, and using a fully rescheduling strategy consumes a lot of computing resources without meeting practical needs. For example, the computational resources consumed by the DPEIA-full algorithm are almost 1.5 times that of the DPEIA algorithm where $\sigma_{q}=0.09$ and qualified 300 products should be completed. Meanwhile, in all experimental situations, the DPEIA algorithm has significant advantages compared to the other four algorithms, validating the competitiveness of the proposed algorithm. Furthermore, in all experimental situations, the DPEIA algorithm consistently outperforms the other four algorithms, validating its competitiveness. The DPEIA algorithm is superior in approximating the Pareto optimality and addressing joint optimization problems with uncertain factors such as quality, production scheduling, and reliability under the QRP-co-effect.

\spacingset{0.9}
\begin{table}[t]
\caption{{Comparison of algorithms performance based on the IGD and HV metrics.}}\label{tablerandom5}
\renewcommand{\arraystretch}{0.9} 
\setlength{\tabcolsep}{2.9mm}{
\begin{tabular}{ccccccccc}
\toprule
Criteria  &Job quality std  &{$N$}   &{EMA} & {MOEA-D} & {MVO} & {NSGAII} & {DPEIA} \\ 
\midrule
   \multirow{6}{*}{IGD}
   &\multirow{3}{*}{  0.03}  &{100} 
  & 1.31 & 0.88 & 0.77 & 0.98 & 0.00 \\
  & &{200} & 0.79 & 0.77 & 0.71 & 1.01 & 0.01  \\
  & &{300} & 0.83 & 0.38 & 0.51 & 0.49 & 0.10  \\
  \cmidrule(r){2-8} 
  &\multirow{3}{*}{  0.09} &{100} & 0.89 & 0.85 & 0.96 & 0.88 & 0.05 \\
  & &{200} & 1.10 & 0.92 & 0.83 & 0.96 & 0.02  \\
  & &{300} & 0.98 & 0.87 & 0.90 & 0.92 & 0.02   \\
  \midrule
  \multirow{6}{*}{HV}
  &\multirow{3}{*}{  0.03} &{100} & 0.00 & 0.07 & 0.08 & 0.03 & 0.41\\
  & &{200} & 0.07 & 0.14 & 0.07 & 0.04 & 0.33 \\
  & &{300} & 0.01 & 0.10 & 0.06 & 0.06 & 0.28 \\
\cmidrule(r){2-8} 
  &\multirow{3}{*}{  0.09}  &{100} & 0.06 & 0.10 & 0.07 & 0.06 & 0.42 \\
  & &{200} & 0.02 & 0.05 & 0.04 & 0.02 & 0.33  \\
  & &{300} & 0.04 & 0.07 & 0.04 & 0.03 & 0.38  \\
\bottomrule
\end{tabular}}
\end{table}
\spacingset{1.5}

  

\spacingset{0.9}
\begin{table*}[t]
\caption{Comparison of algorithms performance based on the RPD-metric.}\label{tablerandom7}
\renewcommand{\arraystretch}{0.85} 
\setlength{\tabcolsep}{7.5mm}{
\begin{tabular}{cccccccccc}
\toprule
 \multirow{2}{*}{Job std ($\sigma_{q}$)} & \multirow{2}{*}{$N$} & \multicolumn{2}{c}{DPEIA-full} & \multicolumn{2}{c}{DPEIA}\\ 
  \cmidrule(r){3-4}\cmidrule(r){5-6} 
   && $C_{max}$ & $C^{m}$ & $C_{max}$ & $C^{m}$\\
\midrule
\multirow{3}{*}{  0.03} & {100} & 4.57 & 0.00 & 0.00 & 0.00\\
& {200} & 2.27 & 8.67  & 2.48 & 15.32\\
& {300}
& 8.41 & 36.60 & 3.92 & 40.23\\
\midrule
\multirow{3}{*}{  0.06} & {100} & 11.16 & 59.97 & 21.19 & 45.11\\
& {200}  & 7.03 & 10.30 & 5.82 & 6.23\\
& {300} & 1.37 & 1.94  & 6.08 & 7.87\\
\midrule
\multirow{3}{*}{  0.09} & {100}  & 14.28 & 32.08 &13.07 & 36.59\\
& {200} & 3.34 & 6.90 & 5.56 & 5.33 \\
& {300} &20.32 & 8.34 & 20.04 & 8.91\\
\bottomrule
\label{tablesen}
\end{tabular}}
\end{table*}
\spacingset{1.5}
\spacingset{0.9}
\begin{table}[t]
\caption{{Comparison of algorithms under different rescheduling schemes and maintenance strategies.}}
\renewcommand{\arraystretch}{0.9} 
\setlength{\tabcolsep}{2.95mm}{
\begin{tabular}{ccccccccccc}
\toprule
{Assessment criteria} &{Job sizes} &{$N$}  & {DPEIA-full} & {DPEIA-single} & {DPEIA}\\
\midrule
  \multirow{3}{*}{  IGD} & {small-sized} &{100}  & 0.21 & 0.48 & 0.20\\
  & {medium-sized} &{200}  & 0.45 & 0.67 & 0.40 \\
  & {large-sized} &{300}  & 0.28 & 0.78 & 0.26 \\
\midrule
  \multirow{3}{*}{  HV} & {small-sized} &{100}  & 0.10 & 0.05 & 0.11\\
  & {medium-sized} &{200}   & 0.06 & 0.03  & 0.06  \\
  & {large-sized} &{300}   & 0.08 & 0.01  & 0.07  \\  
\bottomrule
\end{tabular}}
\label{tablerandom8}
\end{table}
\spacingset{1.5}

\section{Conclusion}
This study aims to address a joint optimization problem that encompasses quality, reliability, and production (QRP) for a multi-component production system incorporating reworking activities. The primary objective of the proposed methodology is to attain high-quality outputs, heightened productivity, and increased production system’s reliability, while simultaneously minimizing costs. The study explores the interdependencies among quality, production scheduling, and machine reliability by considering the uncertainties associated with machine degradation, job processing time, and product quality throughout the production process. To address this problem, leveraging the available conditional information on the platform becomes pivotal for enhancing the efficacy of joint optimization. In light of this, we introduce a dual-module solution framework. Firstly, we harness the dynamic interplay among communication mechanisms, maintenance strategies, and relevant properties to enhance computational efficiency. Secondly, we integrate online data effectively, employing dynamic decision-making processes to derive adaptive outcomes. Ultimately, empirical validation through numerical studies substantiates the efficacy of this approach.

\bibliographystyle{chicago}
\spacingset{0.8}
\bibliography{cas-refs}

\newpage
\section*{Appendices}
\appendix

\section{}
\label{app1}
\setcounter{equation}{0}
\renewcommand\theequation{A.\arabic{equation}}

The proof of obtaining for $N_{i,k}^{pm}$.

To obtain the value of $N_{i,k}^{pm}$, variable $z_{i,k}^{a}$ is introduced to locate the time point of the most recent corrective maintenance (CM) activity for machine $k$, and variable $z_{i,k}^{b}$ is used to determine whether the machine has or not undergone any CM activity before job $i$ is completed. For $\forall i \in N_k, k \in K$, the formula for obtaining $N_{i,k}^{pm}$ is designed as follows:
\begin{alignat}{2}
&\sum_{w=1}^{l}(\sum_{(i^{(1)},...i^{(w)})\in N_{k}}(\prod _{j=1}^{w}x_{i^{(j)},i^{(j-1)}})z_{i^{(w)}}) \leq Mz_{i,k}^{a},\\
&\sum_{w=1}^{l}(\sum_{(i^{(1)},...i^{(w)})\in N_{k}}(\prod _{j=1}^{w}x_{i^{(j)},i^{(j-1)}})z_{i^{(w)}}) \geq z_{i,k}^{a},\\
&\sum_{w=1}^{l}(\sum_{(i^{(1)},...i^{(w)})\in N_{k}}(\prod _{j=1}^{w}x_{i^{(j)},i^{(j-1)}})\sum_{g\in G}y_{i^{(w)},g})= M(1-z_{i^{(l)}}+z_{i,k}^{a})+N_{i,k}^{pn},\\
&y_{i,k}^{e}=\left\{\begin{matrix}
0,\ N_{i,k}^{pn}\leq 0,\\ 
1,\ N_{i,k}^{pn}>0
\end{matrix}\right.,\\
&z_{i,k}^{b}=\left\{\begin{matrix}
0,\ \sum_{w=1}^{l^{0}}(\sum_{(i^{(1)},...i^{(w)})\in N_{k}}(\prod _{j=1}^{w}x_{i^{(j)},i^{(j-1)}})z_{i^{(w)}})=0\\ 
1,\ \sum_{w=1}^{l^{0}}(\sum_{(i^{(1)},...i^{(w)})\in N_{k}}(\prod _{j=1}^{w}x_{i^{(j)},i^{(j-1)}})z_{i^{(w)}})>0
\end{matrix}\right.,\\
&N_{k,j}^{pm}= \sum_{j^{'}=1}^{j}y_{k,j^{'},j}^{e}N_{i,k}^{pn}z_{k,j}^{b}+\sum_{w=1}^{l^{0}}(\sum_{(i^{(1)},...i^{(w)})\in N_{k}}(\prod _{j=1}^{w}x_{i^{(j)},i^{(j-1)}})\sum_{g\in G}y_{i^{(w)},g})(1-z_{k,j}^{b}),\\
&z_{i,k}^{a}\in \{0,1\},
\end{alignat}
where $M=n$, and subscript $i^{(j)}$ represents job $i$ processed on the $j^{th}$ position on the machine. $N_{i,k}^{pn}$ denotes the intermediate variable to obtain $N_{k,j}^{pm}$ when at least one CM activity has been performed until complete job $i$. Index $l$ denotes the position of job $i$ processed on the machine, and index $w$ is used to traverse the positions between the first position and position $l$ of the machine.

\section{}
\setcounter{equation}{0}
\renewcommand\theequation{B.\arabic{equation}}
The proof of \emph{Proposition 1}.

A solution can be evaluated from four perspectives based on the objective function: the completion time point and the qualities of jobs, the degradation state of the machines and the maintenance cost in the process. For solution $So$ where job $i$ and job $i^{'}$ is scheduled on the position $j$ and $j+1$ of machine $k$, respectively, we assume that  $T_{k}^{So}$ is the time required by machine $k$ from the initial time to the completion time point of job $i^{'}$. The alternative solution denotes $So^{'}$ where the sequence of jobs $i$ and $i^{'}$ is reversed taking time $T_{k}^{So^{'}}$ until finishing job $i$. $Q_{k}^{So}$ represents the number of the qualified products. $W_{k}^{So}$ denotes the degradation of machine $k$ after processing job $i^{'}$ in solution $So$, and $W_{k}^{So^{'}}$ denotes the degradation state after processing job $i$ in solution $So^{'}$. Let $W_{i^{''},k}$ represents the machine $k$'s degradation state before processing job $i$ and job $i^{'}$. Meanwhile, all jobs are available in the initial time. Thus, we discuss in two conditions, and \emph{Proposition 1} can be derived as follows:

1) The first condition.

According to equation (\ref{C8}), the input jobs whose initial quality are eligible when the output products are qualified. Then, the degradation formulation (equation (\ref{C6})) consists of the Gaussian distribution (according to the proof of Lemma \ref{Lemma 1}) and the Gamma distribution. Under the additivity and the non-negativity for the gamma and Gaussian distributions along with the independent increment property, $h(\alpha_{k}, \Delta t_{i,k}, \beta_{k}, \sigma_{k})$ denotes the increments in $\Delta t_{i,k}$. Meanwhile, the job's available time can be expressed as $t_{i^{'},k}^{a}=t_{i,k}^{a}=0$.
If $O_{i,k} < O_{i^{'},k}$:
\begin{flalign}
&T_{k}^{So}-T_{k}^{So{'}}\\
&=t_{i,k}^{a}+p_{i}^{So}+p_{i^{'}}^{So}-(t_{i^{'},k}^{a}+p_{i^{'}}^{So^{'}}+p_{i}^{So^{'}})=p_{i}^{So}+p_{i^{'}}^{So}-(p_{i^{'}}^{So^{'}}+p_{i}^{So^{'}})<0,\\
&W_{k}^{So}-W_{k}^{So^{'}}\\
&=(W_{i^{''},k}+(y(p_{i}^{So})+y(p_{i^{'}}^{So})+g(\alpha_{k}, (T_{k}^{So}), \beta_{k}, \sigma_{k})))-(W_{i^{''},k}+(y(p_{i}^{So^{'}})+y(p_{i^{'}}^{So^{'}})))+\notag \\&g(\alpha_{k}, (T_{k}^{So^{'}}), \beta_{k}, \sigma_{k})))\\
&=g(\alpha_{k}, p_{i^{'}}^{So}, \beta_{k}, \sigma_{k})+y(p_{i^{'}}^{So})-(g(\alpha_{k}, p_{i}^{So^{'}}, \beta_{k}, \sigma_{k})+y(p_{i}^{So^{'}}))<0.
\end{flalign}

According to equation (\ref{C8}), the quality characteristic can be obtained as follows: 
\begin{flalign}
&D_{i}^{So}-D_{i}^{So^{'}}\\
&=a_{k}W_{i,k}^{So}+W_{i,k}^{So}\times \Gamma_{k}\times \varepsilon_{i}-(a_{k}W_{i,k}^{So^{'}}+W_{i,k}^{So^{'}}\times \Gamma_{k}\times \varepsilon_{i})\\
&<(W_{i,k}^{So}-W_{i,k}^{So^{'}})(a_{k}+\Gamma_{k}\times \varepsilon_{i})<0,
\end{flalign}
where $p_{i}^{So}=O_{i,k}(1+\eta(W_{i^{''},k}+g(\alpha_{k},  t_{i,k}^{a}, \beta_{k}, \sigma_{k})))$,\\$p_{i^{'}}^{So}=O_{i^{'},k}(1+\eta(W_{i^{''},k}+g(\alpha_{k}, (t_{i,k}^{a}+p_{i}^{So}), \beta_{k}, \sigma_{k})+y(p_{i}^{So}))$, $p_{i^{'}}^{So^{'}}=O_{i^{'},k}(1+\eta(W_{i^{''},k}+g(\alpha_{k},  t_{i^{'},k}^{a}, \beta_{k}, \sigma_{k})))$,\\$p_{i}^{So^{'}}=O_{i,k}(1+\eta(W_{i^{''},k}+g(\alpha_{k}, (t_{i^{'},k}+p_{i^{'}}^{So^{'}}), \beta_{k}, \sigma_{k})+y(p_{i^{'}}^{So^{'}}))$.  $y(p_{So,i})$ denotes normalization equation of $p_{So,i}$, and $y$ is a linear function.

2) The second condition.

Similarly, jobs' effects on the machine's degradation follow the Gaussian distribution. Then, in this situation, the following formulas are also satisfied:
\begin{flalign}
&T_{k}^{So}-T_{k}^{So{'}}<0,\quad W_{k}^{So}-W_{k}^{So^{'}},\\
&\Delta W_{i,k}(a_{k}+\Gamma_{k}\times \varepsilon_{i})<\Delta D_{i^{'}}, \\
&\Rightarrow (W_{i,k}^{So}-W_{i,k}^{So^{'}})(a_{k}+\Gamma_{k}\times \varepsilon_{i}) <\Delta D_{i^{'}},  \\
&\Rightarrow D_{i^{'}}^{So}-D_{i^{'}}^{So^{'}}<\Delta D_{i^{'}}, \\
&\Rightarrow |D_{i^{'}}^{So}-SL_{i^{'}}^{(+)}|<\xi_{i^{'}}\\
&\Rightarrow Q_{k}^{So}=Q_{k}^{So^{'}}=2.
\end{flalign}

Meanwhile, for two output products, the number of conforming products can be expressed as $Q_{i,k,i^{'}}^{So^{'}}=0$ under the second condition, while $0\leq Q_{i,k,i^{'}}^{So}\leq 2$, then $Q_{k}^{So}\geq Q_{i,k,i^{'}}^{So^{'}}$.

\emph{Proposition 1} is proved.
\hfill $\blacksquare$

\section{}
\setcounter{equation}{0}
\renewcommand\theequation{C.\arabic{equation}}
The proof of \emph{Lemma 1}.

According to \cite{karimi2012efficient}, in the parallel system, the active schedule contains an optimum solution where the objective function is makespan. Then the completion time of jobs is shorter when applying for the active schedule. Meanwhile, the longer jobs' completion time can accelerate the machine deterioration and quality of output products jobs based on the proof of \emph{Proposition 1}, where other conditions are kept the same. 

\emph{Lemma 1} is proved.
\hfill $\blacksquare$

\section{}
\setcounter{equation}{0}
\renewcommand\theequation{D.\arabic{equation}}
The diagram of the idle decoding.

\begin{figure*}[t]
\centering
\includegraphics[scale=0.48]{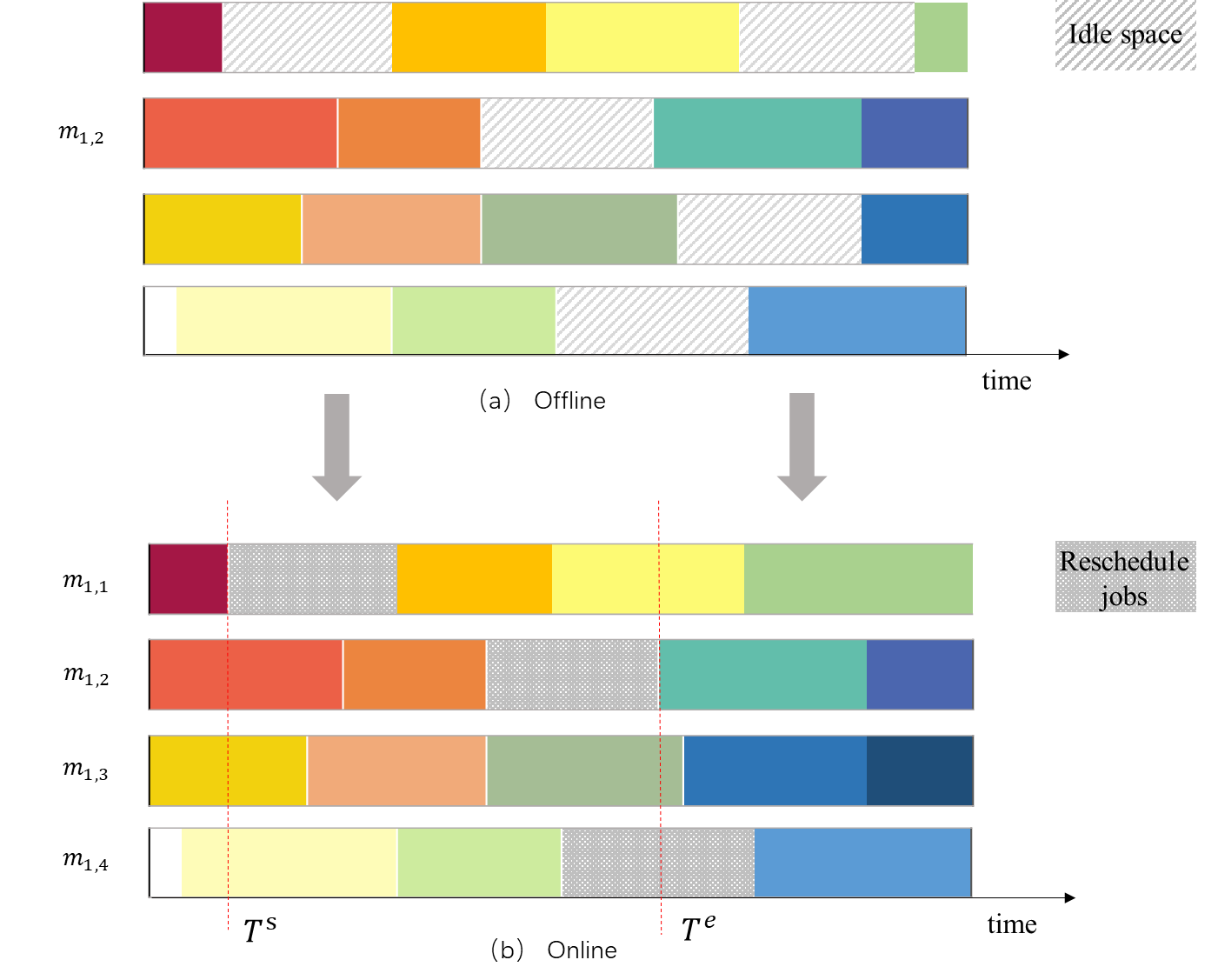}
\caption{\textcolor{black}{The diagram of the idle decoding.}}
\label{Figure_decode}
\end{figure*}
Figure \ref{Figure_decode}(a) shows a schematic Gantt chart of idle spaces, represented by diagonal-stripe rectangles, added during the planning module. In Figure \ref{Figure_decode}(b), plaid rectangles represent rework jobs, while other spaces in different colors denote jobs assigned based on the planned baseline. Jobs between the starting time point $T^{(+)}$ and the ending time point $T^{e}$ are rescheduled. In the evaluation module, it can be observed that the positions of the idle spaces are occupied by unprocessed jobs, as depicted in Figure \ref{Figure_decode}(b).
\newpage

\section{}
\setcounter{equation}{0}
\renewcommand\theequation{E.\arabic{equation}}
The proof of \emph{Proposition 2}.
\begin{flalign}
&\quad \frac{n^{pm}}{n^{c}}<\varsigma^{m},\\
& \quad \Rightarrow (\frac{n^{pm}}{n^{c}}+1)^{2}<(\varsigma^{m}+1)^{2},\\
& \quad \Rightarrow (\frac{n^{c}+n^{pm}}{n^{c}})^{2}<1+2\varsigma^{m}+(\varsigma^{m})^{2},\\
& \quad \Rightarrow (\frac{n^{c}+n^{pm}}{n^{c}})^{2}<1+\frac{C^{pm}}{C^{c}}+\frac{T^{pm}}{T^{c}}+\frac{T^{pm}C^{pm}}{T^{c}C^{c}},\\
& \quad \Rightarrow (\frac{n^{c}+n^{pm}}{n^{c}})^{2}<\frac{T^{c}C^{c}+(T^{c}C^{pm}+T^{pm}C^{c})+T^{pm}C^{pm}}{T^{c}C^{c}},\\
& \quad \Rightarrow (\frac{n^{c}+n^{pm}}{n^{c}})^{2}<\frac{(T^{c}+T^{pm})(C^{c}+C^{pm})}{T^{c}C^{c}},\\
& \quad \Rightarrow \frac{(n^{c})^{2}}{T^{c}C^{c}}>\frac{(n^{c}+n^{pm})^{2}}{(T^{c}+T^{pm})(C^{c}+C^{pm})}.
\end{flalign}

\emph{Proposition 2} is proved.
\hfill $\blacksquare$

\section{}
\setcounter{equation}{0}
\renewcommand\theequation{D1.\arabic{equation}}
\begin{figure*}[h]
\centering
\includegraphics[scale=0.55]{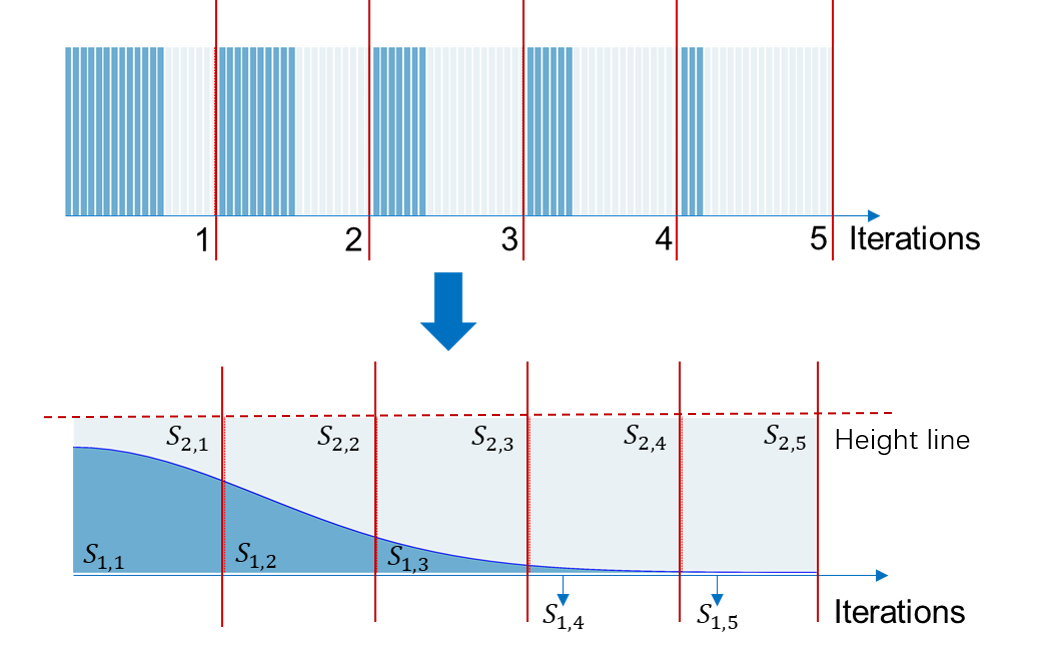}
\caption{\textcolor{black}{The diagram of the computing resource allocation.}}
\label{percent}
\end{figure*}
\hfill $\blacksquare$

\section{}
\setcounter{equation}{0}
\renewcommand\theequation{D0.\arabic{equation}}
Procedure of the EMODE algorithm.
\spacingset{1.0}
\begin{algorithm}[]
  \caption{Procedure of the EMODE algorithm}
  \label{Algorithm 1}
  \LinesNumbered
  \KwIn{$N$: the set of job indexes, \\$K$: the set of machine indexes,\\ $popSize^{P}$: the population size in the planning module,\\ $maxIter^{P}$: the number of iterations in the planning module
}
  \KwOut{Planning}

  \textbf{Initialize:} $g=0$, $popInit$.initialize()\\
    
\baselineskip=0pt

  \While{$g<maxIter^{P}$}{
       Select the number $\nu$ from 2 to 0 based on formula (\ref{a_1});\\
      \If{$\nu>1$}
      {
            \For{$i = 1$ to $popSize^{P}$}
            { 
            \If {$rand < 0.7$}
                {Use the similarity-oriented RE operator to generate a new solution;}
            \Else
                {Use the decoding-based DE operator to generate a new solution;}
            }
        }
        
      \Else{
            \For{$i = 1$ to $popSize^{P}$}
            {
                randomly select two jobs assigned on two consecutive positions $j$ and $j^{'}$ in a maintenance interval belonging to a machine;\\ 
                \If {\textnormal{these two jobs satisfy one of sub-conditions in \emph{Proposition 1}}}
                {
                    generate new individual $x^{'}$ by swapping two jobs' positions, and $popInit \cup x^{'}$;
                }
                \Else
                {
                calculate the ratio of total processing time on each machine and makespan to obtain the busiest machine ${m^{b}}$ and the idlest machine ${m^{d}}$;\\
                    \If{\textnormal{${m^{b}}$ is not empty}}
                    {
                        generate new individual $x^{'}$ by randomly selecting a job on the most busiest machine ${m^{d}}$ to assign to the idlest machine ${m^{b}}$; \\
                        $popInit \cup x^{'}$;\\
                    }
                }
                    }
            
            }

        update population $popInit$ based on the roulette method with fitness function(\ref{fit_1}) value as probability and the $popSize^{P}$;\\
        
        $g$++;\\
   }
\end{algorithm} 
\hfill $\blacksquare$
\spacingset{1.5}

\section{}
\setcounter{equation}{0}
\renewcommand\theequation{D5.\arabic{equation}}
Procedure of the EMODE rescheduling.
\spacingset{1.0}
\begin{algorithm}[]
    \baselineskip=0pt
  \caption{Procedure of rescheduling}
  \label{Algorithm 2}
  \LinesNumbered
  \KwIn{$N$: the set of job indexes, \\$K$: the set of machine indexes,\\ $popSize^{R}$:the population size in the online improvement module, \\
  $baseline$: the current schedule plan,\\ $maxIter^{R}$: the number of iterations in the online improvement module }
  \KwOut{$Schedule$}
  \textbf{Initialize:} $g=0$, $popRes$.initialize()\\
  
  \baselineskip=0pt
  
  \While{$g<maxIter^{R}$}
  {     \If {$g$ mod 2==1}
        {
            determine machine $m^{key}$ effecting the makspan and randomly select a job $i$ assigned on the machine;\\
            determine machine $m^{erl}$ with the earliest completion time;\\
            rearrange job $i$ to a position on machine $m^{erl}$ satisfying the following condition;\\
            \For{$k = 1$ to $n^{inr}_{m^{erl}}$ - 1} 
            {
                \For{$j$ = 1 to $n^{inr}_{m^{erl},u}$ - 1}
                {
                    \If {$p_{i}$ $>=$ $p_{i^{m^{erl}_{j}}}$ and $p_{i}$ $<=$ $p_{i^{m^{erl}_{j+1}}}$}
                    {
                        insert the job $i$ into the position $j$ on machine $m^{erl}$, and generate new individual $x^{'}$;\\
                        $popRes \cup x^{'}$;\\
                        pass;\\
                    }
                }
            }
        }
        \Else
        {
            \For {$m$ in $K$}
            {
                randomly select positions $j$ and $j^{'}$ on machine $m$ and swap the jobs assigned on these positions, and generate new individual $x^{'}$;\\
                $popRes \cup x^{'}$;\\
            }
        }
            randomly select two positions on different machines and swap the processing sequence of these two jobs;\\
        $g++$;
  }
  Notes: $n^{inr}_{m^{erl}}$ and $n^{inr}_{m^{erl},u}$ represent the number of intervals divided by maintenance activities and the $u^{th}$ interval for machine $m^{erl}$, respectively; $p_{i^{m^{erl}_{j}}}$ denotes the actual processing time of the job processed on $j^{th}$ position in the $u^{th}$ interval of machine $m^{erl}$.
\end{algorithm}
\hfill $\blacksquare$
\spacingset{1.5}
\section{}
\setcounter{equation}{0}
\renewcommand\theequation{D2.\arabic{equation}}
\begin{table*}[h]
\caption {\label{tablepa_3} {Parameter setups of each machine for the base case.}}
\setlength{\tabcolsep}{6.92mm}{
\begin{tabular}{ccccc}
\toprule
 Parameters & $M_{1}$ & $M_{2}$ & $M_{3}$ & $M_{4}$\\ 
\midrule
  $\mu_{k}^{(-)}$ & 82.4 & 66.4 & 74.72 & 66 \\
  $\sigma_{k}^{(-)}$  & 0.00306 & 0.00296 & 0.00326 & 0.00254  \\
  $\mu^{(+)}_{k}$ & 0.0 & 0.0 & 0.0 & 0.0  \\
  $\sigma^{(+)}_{k}$  & 0.015 & 0.015 & 0.015 & 0.015  \\
  $a_k$  & 91.1 & 98.95 & 103.5 & 86.5\\
  $b_k$  & 0.57032 & 0.5664 & 0.5832 & 0.5612\\
  $\beta_k$ & 5.792e-05 & 5.516e-05 & 6.423e-05 & 6.085e-05\\
  $C^{pm}_{k}$ & 430 & 275 & 230 & 195\\
  $C^{cm}_{k}$ & 1312 & 1028 & 876 & 832\\
  $W^{0}_{k}$ & 0.1 & 0.105 & 0.11 & 0.99\\
  $L_k$ & 0.35 & 0.4025 & 0.385 & 0.315\\
  $T^{ps}_{k}$ & 12.6 & 10.85 & 10.5 & 10.15\\  
  $T^{pm}_{k}$ & 12.54 & 10.92 & 10.49 & 10.15\\
  $T^{cm}_{k}$ & 44.75 & 40.50 & 36.64 & 36.64\\
\bottomrule
\end{tabular}}
\end{table*}

\begin{center}
$\boldsymbol a^{'}= \left[ 
\begin{array}{cccccc}
	& 0.0112 & 0.0173 & 0.0147 & 0.0158& \\
\end{array}
\right],$\\

$\boldsymbol b^{0'}= \left[ 
\begin{array}{cccccc}
	& 0.0098 & 0.0106 & 0.0105 & 0.0072& \\
\end{array}
\right],$\\

$\Gamma^{'}= \left[ 
\begin{array}{cccccc}
	& 0.0137 & 0.0152 & 0.0132 & 0.0143& \\
\end{array}
\right],$\\


$\mu^{'}= \left[ 
\begin{array}{cccccc}
	& 42.72 & 42.72 & 42.72 & 42.72 \\
\end{array}
\right],$\\

$SL^{(+)}_{e^{{th}^{'}}}= \left[ 
\begin{array}{cccc}
	& 42.72 & 42.61&\\
\end{array}
\right],$ \quad
$\xi_{e^{th}}^{'}= \left[ 
\begin{array}{cccc}
	& 0.08 & 0.07&\\
\end{array}
\right].$\\
\end{center}
\hfill $\blacksquare$

\section{}
\setcounter{equation}{0}
\renewcommand\theequation{D3.\arabic{equation}}
The DPEIA algorithm is equipped with the following parameter settings: the mean $\mu_{c}$ and qualified deviation $\sigma_{c}$ of the truncated normal distribution, the communication mechanism control parameter $\varpi$ and the degradation coefficient are set to 0, 1.13, and 0.5, respectively. The degradation coefficient $\eta$ is set to 0.2.
The maximum iteration number $maxIter^{P}$ for all algorithms in the experiment is set to 100. Given the stochastic nature of the problem, each solution is repeated 50 times to obtain average results. Furthermore, the imperfect maintenance strategy incorporates influence coefficients ($\theta$ and $\varphi$) set to 0.2 and 0.08, respectively. 

1) IGD-metric:
\begin{flalign}
IGD\left(P, P^*\right)=\frac{1}{\left|P^*\right|} \sum_{p \in P^*} \operatorname{dist}(p, P), 
\end{flalign}
where $P$ and $P^*$ represent the set of points on the Pareto optimal surface obtained by the corresponding algorithm and all algorithms, respectively, and $p$ = $(obj_1^*, obj_2^*)$. Value $obj_i^*$ is the normalized form of objective function values $obj_i$, namely, the makespan and the cost in this study: 
\begin{flalign}
& obj_i^*=\frac{Obj_i-Obj_{i, min }}{Obj_{i, max }-Obj_{i, min }},
\end{flalign}
where $Obj_{i, min }$ and $Obj_{i, max }$ are the minimum and maximum values of objective $obj_i$ in total algorithms.

2) hypervolume-metric:
\begin{flalign}
HV=\delta\left(\bigcup_{i=1}^{|S|} r_i\right),
\end{flalign}
where $S$ denotes the set of non-dominated solutions obtained by the corresponding algorithm, while $r_i$ represents the area of the rectangle constructed using the normalized objective space in the range of $[0,1]$ and the qualified point $(1,1)$ as the opposite vertices, which represents the worst solution among all algorithms. Symbol $\delta$ measures the volume of the contents in parentheses. Furthermore, a larger value of the hypervolume indicator ($HV$) corresponds to a better performance of the algorithm.

3) RPD-metric:
\begin{flalign}
RPD_i= 100 \times \frac{Obj_{i,method} - Obj_{i,best}}{Obj_{i,best}},
\end{flalign}
where $Obj_{i,method}$ represents the average value of the non-dominated solutions attained by the algorithm in question on the $i^{th}$ objective. $Obj_{i,best}$ refers to the best value of the non-dominated solutions found by all the algorithms on the same $i^{th}$ objective. 
\hfill $\blacksquare$


\section{}
\setcounter{equation}{0}
\renewcommand\theequation{F.\arabic{equation}}
The encoding and the decoding parts of DPEIA.

\begin{figure}[h]
\centering
\includegraphics[scale=0.6]{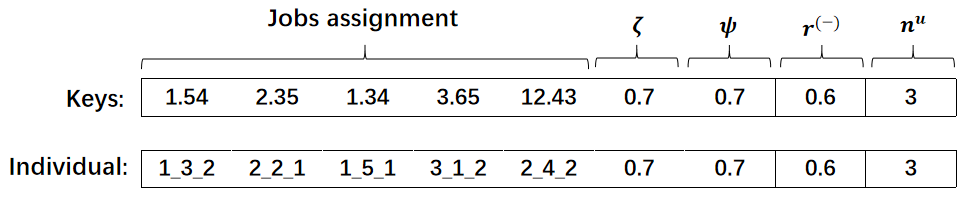}
\caption{\textcolor{black}{An illustration using keys and a chromosome in the random-key EMODE.}}
\label{Figure_Gachorm}
\end{figure}

\noindent {\bf\emph{ Step 1. Encoding.}} All jobs are allocated to machines, forming an individual with five components. The initial part encompasses job assignments to machines, including the job sequence. The second component represents the threshold of the Preventive Maintenance (PM) coefficient $\zeta$. The third part indicates the PM group proportion associated with the coefficient $\psi$. The fourth part involves the rescheduling triggering threshold linked to the coefficient $Thr^{r}$. The fifth part denotes the maximum number of times a machine's PM coefficient $n^{u}$.
In the primary component concerning job assignments, the integer part designates the assigned machine, while the decimal part determines the job order, similar to idle time. An illustrative example of an individual is presented in Figure \ref{Figure_Gachorm}. In the section related to job assignments, the integer part of the numbers signifies the number of assigned machines (with idle time being the remainder after division by 10). The decimal parts are sorted in ascending order, and the sequence is employed as the processing sequence of jobs on the machine. In the chromosome section, $a_{b_c}$ signifies that job $b$ is allocated to the $c^{th}$ position on machine $a$.

\noindent {\bf\emph{ Step 2. Decoding.}} Utilizing the established maintenance policy and the designed individual, we employ a Monte Carlo simulation to replicate the system's stochastic deterioration process. This simulation is instrumental in obtaining the total cost and makespan based on the decision variables $\zeta$, $\psi$, $Thr^{r}$, and $n^{u}$. The algorithm's detailed procedure is outlined in Algorithm \ref{Algorithm 1}, where $N^{s}$ represents the sample size of the simulation.\\

\noindent {\bf\emph{Step 3. Parameters Updating.}} The fundamental steps of EMODE, encompassing individual selection, reproduction, crossover, and mutation, are executed to update the job assignment, sequencing plan, and the maintenance decision variables $\zeta$, $\psi$, $Thr^{r}$, and $n^{u}$.\\

\noindent {\bf\emph{Step 4. Solution.}} Steps 2 and 3 are iterated until the termination criterion is satisfied. Subsequently, the optimal job assignment and sequencing plan can be determined. Additionally, the optimal decision variables for OM, i.e., $\zeta$, $\psi$, $Thr^{r}$, and $n^{u}$, can be obtained.

\end{document}